\documentclass[10pt,journal]{IEEEtran}
\usepackage[utf8]{inputenc}
\usepackage{amsmath,amssymb,amsfonts}
\usepackage{graphicx}
\usepackage{color}
\usepackage{algorithmic}
\usepackage{textcomp}
\usepackage{bm}
\usepackage{epstopdf}  
\usepackage[mathscr]{eucal}
\usepackage{enumerate}
\usepackage{cite}

\usepackage[colorlinks=true,  linkcolor=black, citecolor=red, breaklinks=true, urlcolor=blue]{hyperref}

\graphicspath{{./figures/includeTyler/}{./figures/}{./figures/tikz/}}
\usepackage{tikz}
\usepackage{pgfplots}
\pgfplotsset{compat=1.13}
\newlength\fwidth

\newcommand{\be}{\begin{equation}}
\newcommand{\ee}{\end{equation}}
\renewcommand{\Vec}[1]{\boldsymbol{#1}}
\newcommand{\Mat}[1]{\boldsymbol{#1}}
\newcommand{\Exp}[1]{\mathrm{e}^{#1}}
\newcommand{\trace}[1]{\mathop{\mathrm{tr}}\left[{#1}\right]}

\newcommand{\delete}[1]{}

\newcommand{\abstractbegin}{\textcolor{red}{\bf TSP draft \today\ }}

\newcommand{\cfmcomment}[1]{\textcolor{blue}{(comment by cfm: {#1})}}
\newcommand{\esacomment}[1]{\textcolor{magenta}{(comment by Esa: {#1})}}

\renewcommand{\abstractbegin}{}

\renewcommand{\cfmcomment}[1]{}
\renewcommand{\esacomment}[1]{}

\begin{document}

\title{Robust and Sparse M-Estimation of DOA}

\author{Christoph F. Mecklenbr\"auker,~\IEEEmembership{Senior Member,~IEEE,} Peter Gerstoft,~\IEEEmembership{Fellow,~IEEE,}
Esa Ollila~\IEEEmembership{Senior Member,~IEEE},
Yongsung Park~\IEEEmembership{Member,~IEEE }%
\thanks{Manuscript received Dec.  22,  2022; revised May31,  2023. This work was supported by Office of Naval Research.}%
\thanks{C.F. Mecklenbr\"auker is with Inst. of Telecommunications, TU Wien, Vienna, Austria.} 
\thanks{Y. Park and P. Gerstoft are with NoiseLab, UCSD, San Diego (CA), USA.}
\thanks{E. Ollila is with Dept.  of Information and Communications Engineering,  Aalto University,  Aalto,  Finland.}
\thanks{A preliminary version of the proposed DOA M-estimation algorithm and obtained results were published in \cite{Mecklenbraeuker2022icassp}.
Codes  are available at GitHub \cite{RobustSBL-github}.}
}
\maketitle
\begin{abstract}
\abstractbegin
A robust and sparse Direction of Arrival (DOA) estimator is derived  for array data that follows a Complex Elliptically Symmetric (CES) distribution with zero-mean and finite second-order moments.
The derivation allows to choose the loss function and four loss functions are discussed in detail: the Gauss loss which is the Maximum-Likelihood (ML) loss for the circularly symmetric complex Gaussian distribution,  the ML-loss for the complex multivariate $t$-distribution (MVT) with $\nu$ degrees of freedom,  as well as  Huber  and Tyler loss functions.   For Gauss loss,  the method reduces to Sparse Bayesian Learning (SBL).
The root mean square DOA error of the derived estimators is discussed for Gaussian,  MVT,  and $\epsilon$-contaminated data.
The robust SBL estimators perform well for all cases and nearly identical with classical SBL for Gaussian noise.
\end{abstract}



\section{Introduction}
\label{sec:intro}

Heavy-tailed sensor array data arises e.g.  due to clutter in radar \cite{Gini-Greco2002} and interference in wireless links  \cite{Clavier2021}. 
Such array data demand statistically robust array processing. 
There is a rich literature on statistical robustness
\cite{Huber2011,Hampel2011,Zoubir-Koivunen-Chakhchoukh-Muma2012,Zoubir2018}.  
In this work,  we derive,  formulate,  and investigate a robust statistical method that emulates Sparse Bayesian Learning (SBL) for DOA from array data \cite{gerstoft2016mmv},   but is not unduly affected by outliers or other small departures from data distribution assumptions.  The proposed method assumes Complex Elliptically Symmetric (CES) data.
In engineering applications,  CES data models have been widely used when non-Gaussian models are needed \cite{ollila2012complex,Feintuch2023,Besson2016,Singh2021,Zhang2016sirp,Meriaux2019sirp,Luo2022deep-mimo,viberg1997}. 

Due to the central limit theorem,  data are often modeled as Gaussian.  SBL was derived under a joint complex multivariate Gaussian assumption on source signals and noise \cite{Tipping2001}.  Direction of arrival (DOA) estimation for plane waves using SBL is proposed in Ref.~\cite[Table~I]{gerstoft2016mmv}. 
In some cases,  the data contain outliers and robust processing is desired to ensure that good estimates are computed even when these rare events occur.

SBL provides DOA estimates based on the sample covariance matrix (SCM) of the data sample.
There is a rich literature on DOA estimation based on SBL for single and multiple measurements vectors
\cite{Wipf2007beam,yang2012off,dai2017sparse,DaiSo2018,pote2020robustness,pote2023maximum,wu2015,wang2015,chen2018,Nannuru2019,yang2020,wang2021}:
In  Ref.\ \cite{Wipf2007beam} the SCM  in the minimum variance adaptive beamformer was replaced with an estimate from a relevance vector machine.  This improves its performance in case of source correlations or limited number of snapshots.
An iterative algorithm using a Bayesian approach for an off-grid DOA model with joint sparsity among the snapshots was developed in
Ref.\ \cite{yang2012off}.
SBL-type algorithms for DOA estimation in the presence of impulsive noise were proposed in  Ref.\ \cite{dai2017sparse,DaiSo2018}.
This Bayes-optimal approach has high resolution and accuracy at the cost of high numerical complexity when the DOA grid is dense. 

The SCM is a sufficient statistic under the jointly Gaussian assumption, but it is not robust against deviations from this assumption \cite{Ollila2003pimrc}. 
The SBL approach is flexible through the usage of various priors, although Gaussian are most common \cite{gerstoft2016}.
For Gaussian priors this has been approached based on minimization-majorization \cite{Stoica2012} and with expectation maximization (EM) \cite{Wipf2007,Wipf2007beam,Wipf2004,Zhang2011,Liu2012, Zhang2016}. 
We estimate the hyperparameters iteratively from the likelihood derivatives using stochastic maximum likelihood\cite{Boehme1985,Jaffer1988,Stoica-Nehorai1995}.
A numerically efficient SBL implementation is available on GitHub \cite{SBL4-github}.
Recent investigations showed that Sparse Bayesian Learning (SBL) is 
lacking in statistical robustness \cite{Ollila2003pimrc,Mecklenbraeuker-Gerstoft-Ollila-WSA2021,Mecklenbraeuker2022icassp}. 
A Bayes-optimal algorithm was proposed to estimate DOAs in the presence of impulsive noise from the perspective of SBL in \cite{DaiSo2018}.  In the following, we derive robust and sparse Bayesian learning which can be understood as introducing a data-dependent weighting into the SCM.

Previously,  a DOA estimator for plane waves observed by a sensor array based on a complex multivariate Student $t$-distribution data model was studied.
A qualitatively robust and sparse DOA estimate was derived as Maximum Likelihood (ML) estimate based on this model \cite[Sec.~5.4.2]{Zoubir2018},\cite{Ollila2003pimrc,Mecklenbraeuker-Gerstoft-Ollila-WSA2021}.

Here,  we solve the DOA estimation problem from multiple array data snapshots in the SBL framework 
\cite{Wipf2007,gerstoft2016} and use the maximum-a-posteriori (MAP) estimate for DOA reconstruction. 
We assume a CES data model with unknown source variances for the formulation of the likelihood function.
To determine the unknown parameters,  we maximize a Type-II  likelihood (evidence) for this CES  data model
and estimate the hyperparameters iteratively from the likelihood derivatives using stochastic maximum likelihood.

The major contributions of this paper are:
(1) We formulate an SBL algorithm with a general loss function for DOA M-estimation which, given the number of sources,  estimates the set of DOAs corresponding to non-zero source power.  
The data model parametrization is independent of the number of snapshots.
(2) In addition to conventional Gaussian loss,  we investigate the RMSE performance of DOA M-estimation for loss functions associated with data priors featuring potentially strong outliers (Huber,  MVT and Tyler).
This leads to a robust and sparse DOA M-estimator.
(3) The robust and sparse DOA M-estimator is shown to be insensitive to heavy tails, outliers,  and unknown correlations among the sources.

The outline of the paper is as follows: We introduce the notation and the array data model in Sec.~\ref{sec:model}.
Thereafter, we formulate the objective function and specific loss functions used for DOA M-estimator in Sec.~\ref{sec:M-estimation} and describe the proposed algorithm.  
Simulation results for DOA estimation are discussed in Sec.~\ref{sec:results} and report on 
convergence of the algorithm and associated run time in Sec.~\ref{sec:convergence}.

\section{Complex elliptically symmetric data model}
\label{sec:model}

Narrowband waves are observed on $N$ sensors for $L$ snapshots $\Vec{y}_{\ell} $ and the array data is 
$\Mat{Y}=[\Vec{y}_1\ldots\Vec{y}_L]\in\mathbb{C}^{N\times L}$.
We model the snapshots $\Vec{y}_{\ell}$ by a scale mixture of Gaussian distributions, which have a stochastic decomposition of the form
\begin{align}
\Vec{y}_{\ell} = \sqrt{\tau_{\ell}} \, \Vec{v}_{\ell}, \text{ with }
\Vec{v}_{\ell} = \Mat{A}\Vec{x}_{\ell} + \Vec{n}_{\ell},  \quad(\ell=1\ldots L) \label{eq:scale-mixture}
\end{align}
where $\tau_{\ell}>0$ is a random variable independent of $\Vec{v}_{\ell}$.
The random variables $\tau_{\ell}$ enable to accommodate various array data distributions in the model to be discussed further below.

The data model \eqref{eq:scale-mixture} has been used for modelling heavy-tailed non-Gaussian data, most notably clutter in radar \cite{Feintuch2023,Besson2016,Singh2021,Zhang2016sirp,Meriaux2019sirp,Luo2022deep-mimo,viberg1997}. 
We use CES assumptions differently from 
\cite{Feintuch2023,Besson2016,Singh2021,Zhang2016sirp,Meriaux2019sirp,Luo2022deep-mimo,viberg1997} where the model has been used to model noise, clutter, or interference.
We assume that the snapshots follow a CES distribution.
The motivation for the data model \eqref{eq:scale-mixture} is to formulate a model which includes the additive white Gaussian noise (AWGN) model as a special case and is general enough to include cases with heavy-tailed data and data with  outliers.

A possibility would be to estimate the $\tau_\ell$ for each snapshot.
However,  we  formulate a robust statistical method that is insensitive to outliers and small departures from the model  \eqref{eq:scale-mixture}.

The unknown zero-mean complex source amplitudes are the elements of $\Mat{X}=[\Vec{x}_1\ldots\Vec{x}_L]\in\mathbb{C}^{M\times L}$ where
$M$ is the considered number of hypothetical DOAs on the given grid $\{ \theta_1,\ldots,\theta_M\}$. The source amplitudes are independent across sources and snapshots, i.e.,\ $x_{ml}$ and $x_{m'l'}$ are independent for $(m,l) \ne (m',l')$.
If $K$ sources are present in the $\ell$th array data snapshot, the $\ell$th column of $\Mat{X}$ is  $K$-sparse and we assume that the sparsity pattern is the same for all snapshots.
The sparsity pattern is modeled by the active set 
\begin{align}
\mathcal{M} &=\left\{m\in\{1,\ldots,M\}\,|\,x_{m\ell}\ne0\right\}=\{m_1,\ldots, m_K\},
  \label{eq:mathcal-M}
\end{align}
and $x_{m\ell}=0$ for all $m\not\in\mathcal{M}$.
The noise $\Vec{N}=[\Vec{n}_1\ldots\Vec{n}_L]\in\mathbb{C}^{N\times L}$ is assumed independent identically distributed (iid) across sensors and snapshots, zero-mean, with finite variance $\sigma^2$ for all $n,\ell$. 

The $M$ columns of the dictionary $\Mat{A}=[\Vec{a}_1\ldots\Vec{a}_M]\in\mathbb{C}^{N\times M}$ are the replica vectors for all hypothetical DOAs.  For a uniform linear array (ULA), the dictionary matrix elements are  $A_{nm}=\mathrm{e}^{-\mathrm{j}(n-1)\frac{2\pi d}{\lambda}\sin\theta_m}$ ($d$ is the element spacing and $\lambda$ the wavelength). 
The $K$ ``active'' replica vectors are aggregated in 
\begin{align}
\Mat{A}_{\mathcal{M}} = [\Vec{a}_{m_1} \ldots\Vec{a}_{m_K}] \in\mathbb{C}^{N\times K},
\end{align}
with its $k$th column vector $\Mat{a}_{m_k}$, where $m_k\in\mathcal{M}$.

The source and noise amplitudes are jointly Gaussian and independent of each other, i.e. $\Vec{x}_{\ell} \sim \mathbb{C}\mathcal{N}_M(\Vec{0},\Mat{\Gamma})$ and $\Vec{n}_{\ell} \sim \mathbb{C}\mathcal{N}_N(\Vec{0},\sigma^2\Mat{I}_N)$. It follows from \eqref{eq:scale-mixture}, that  $\Vec{v}_{\ell} \sim \mathbb{C}\mathcal{N}_N(\Vec{0},\Mat{\Sigma})$ with 
\begin{align}
 \Mat{\Sigma} &=   \Mat{A} \Mat{\Gamma} \Mat{A}^{\sf H} + \sigma^2 \Mat{I}_N , \label{eq:Sigma-model_new}\\
\Mat{\Gamma} &=  \mathrm{cov}(\Vec{x}_{\ell}) = \mathop{\mathrm{diag}}(\Vec{\gamma})
\label{eq:Gamma-model_new}
\end{align}
where $\Vec{\gamma} = [\gamma_1 \ldots \gamma_M]^{\sf T}$ is the $K$-sparse vector of unknown source powers.  

The matrix $\Mat{\Sigma}$ can be interpreted as the covariance matrix  $\mathrm{cov}(\Vec{v}_{\ell})$ of the Gaussian component $\Vec{v}_{\ell}$, but this is not observable in this model and the sensor array only observes the scale mixture $\Vec{y}_{\ell}$. The matrix $\Mat{\Sigma}$ is called the scatter matrix in the following.

Since $\Vec{y}_{\ell} | \tau_{\ell} \sim \mathbb{C} \mathcal N_N ( \Vec{0}, \tau_{\ell} \Mat{\Sigma})$, the density of $\Vec{y}_{\ell}$ is 
\begin{align}
p_{\Vec{y}}(\Vec{y}_{\ell}) &= \int_0^\infty \!\!\! p_{\Vec{y},\tau}(\Vec{y}_{\ell},\tau) \,\mathrm{d}\tau = \int_0^\infty \!\!\!  p_{\Vec{y}|\tau}(\Vec{y}_{\ell}|\tau) p_{\tau}(\tau)  \, \mathrm{d}\tau  \\ 
&= (\det\Mat{\Sigma})^{-1} g(\Vec{y}_{\ell}^{\sf H} \Mat{\Sigma}^{-1} \Vec{y}_{\ell}) \label{eq:ces}
 \end{align} 
where the so-called \emph{density generator} $g(\cdot)$ 
\cite{ollila2012complex},\cite[Eq.(4.15)]{Zoubir2018} is evaluated by 
 \be \label{eq:t_g}
 g(t)= \pi^{-N} \int_0^\infty \tau^{-N} e^{-t/\tau} p_{\tau}(\tau)  \,\mathrm{d}\tau.
\ee
The form of \eqref{eq:ces} shows that the distribution of $\Vec{y}_{\ell}$ is CES with mean zero \cite{Richmond,greco2014maximum,Fortunati2019b}. 
The CES distribution generalizes  the complex Gaussian by allowing higher or lower tail probabilities than the complex Gaussian while keeping the elliptical contours  of equal density \cite{ollila2012complex}.

If the random scaling $\sqrt{\tau_{\ell}}=1$ in \eqref{eq:scale-mixture} for all $\ell$ then the commonly assumed Gaussian data model is recovered,  
\be
   \Vec{y}_{\ell} = \Mat{A}\Vec{x}_{\ell} + \Vec{n}_{\ell}~.  \label{eq:linear-model3}
 \ee
This model results in Gaussian data,  $\Vec{y}_{\ell}\sim\mathbb{C}\mathcal{N}_N(\Vec{0},\Mat{\Sigma})$.
Model \eqref{eq:linear-model3} was assumed in \cite{gerstoft2016mmv} while here we assume the more general model \eqref{eq:scale-mixture}.

In array processing applications,  the complex Multi-Variate $t$-distribution (MVT distribution) \cite{kent1991redescending,ollila2020shrinking}  can be used as an alternative to the Gaussian distribution 
in the presence of outliers because it has heavier tails than the Gaussian distribution.  
The MVT-distribution is a suitable choice for such data and provides a parametric approach to robust statistics \cite{Zoubir2018,Mecklenbraeuker-Gerstoft-Ollila-WSA2021}.
The complex MVT distribution is a special case of the CES distribution, for details see Appendix \ref{sec:MVT-array-data-model}.

For numerical performance evaluations of the derived M-estimator of DOA,  three data models are used in Sec.  \ref{sec:results}: Gaussian, MVT, and $\epsilon$-contaminated.  The Gaussian and MVT models are CES, whereas the $\epsilon$-contaminated model is not.

\section{M-estimation based on CES distribution}
\label{sec:M-estimation}

\subsection{Covariance matrix objective function}\label{sec:covmat}
We follow a general approach based on loss functions and assume that the  array data $\Mat{Y}$ has a CES distribution with zero mean $\Vec{0}$ and positive definite Hermitian $N \times N$ covariance matrix parameter $\Mat{\Sigma}$ \cite{ollila2012complex,Mahot2013}. Thus
\begin{align}
p(\Mat{Y}| \Vec{0}, \Vec{\Sigma}) &= \prod_{\ell=1}^L \det(\Vec{\Sigma}^{-1})g(\Vec{y}_{\ell}^{\sf H}\Vec{\Sigma}^{-1}\Vec{y}_{\ell}).
\label{eq:CES-pdf} 
\end{align}
An M-estimator of the covariance matrix $\Mat{\Sigma}$ is defined as a positive definite Hermitian $N \times N$ matrix that minimizes the objective function \cite[(4.20)]{Zoubir2018},
\begin{align}
\mathcal{L}(\Mat{\Sigma})
 &= \frac{1}{L b}\ \sum\limits_{\ell=1}^{L} \rho(\Vec{y}_{\ell}^{\sf H} \Mat{\Sigma}^{-1} \Vec{y}_{\ell})  - \log\det(\Mat{\Sigma}^{-1}) ,
 \label{eq:Mobjective}
\end{align}
where $ \Vec{y}_{\ell} $ is the $\ell$th array snapshot and $\rho: \mathbb{R}_0^+ \to  \mathbb{R}^+$, is called the loss function. The loss function is any continuous, non-decreasing function which satisfies that $ \rho(e^x)$ is convex in $-\infty <x < \infty$, cf. \cite[Sec. 4.3]{Zoubir2018}.  
A specific choice of loss function $\rho$ renders \eqref{eq:Mobjective} equal to the negative log-likelihood of $\Mat{\Sigma}$ when the array data are CES distributed with density generator $g(t)=\mathrm{e}^{-\rho(t)}$ \cite{Huber1964}.  If the loss function  is chosen, e.g.,  as $\rho(t)=t$ then \eqref{eq:Mobjective} becomes the negative log-likelihood function for $\Mat{\Sigma}$ in the Gaussian model \eqref{eq:linear-model3}.

The term $b$ is a fitting coefficient,  called consistency factor,  which 
makes the minimizer of \eqref{eq:Mobjective} consistent to  $\Mat{\Sigma}$ when the array data are Gaussian \cite[Sec. 4.3]{Zoubir2018}.  Thus,  $\Vec{y} \sim \mathbb{C}\mathcal N_N(\Vec{0},\Mat{I})$,  
\begin{align} 
b &= \mathsf{E}[\psi(\| \Vec{y} \|^2)]/N,  \quad\label{eq:consitency factor}  \\ 
&=  \frac 1 N \int_0^\infty \psi( t/2)  f_{\chi^2_{2N}}(t) \mathrm{d} t \label{eq:consitency factor2} 
\end{align} 
 where $\psi(t)=t\, \mathrm{d}\rho(t)/\mathrm{d}t$ and  $f_{\chi^2_{2N}}(t)$ denotes the pdf of chi-squared distribution with $2N$ degrees of freedom. 
To arrive from \eqref{eq:consitency factor} to \eqref{eq:consitency factor2} we used that $\| \Vec{y} \|^2 \sim (1/2)  \chi^2_{2N}$. 

Minimizing \eqref{eq:Mobjective} with $b$ according to \eqref{eq:consitency factor2}  results in a consistent M-estimator of the covariance matrix $\Mat{\Sigma}$ when the objective function is derived under a given non-Gaussian array data assumption (as in Sec. \ref{sec:loss-func}) but is in fact Gaussian ($\Vec{y}_{\ell} \sim \mathbb{C} \mathcal N_N(\Vec{0},\Mat{\Sigma}))$. 
A derivation of the consistency factor $b$ for the loss functions used in this paper is given in Appendix \ref{sec:consistency-factor}.

\subsection{Loss functions}
\label{sec:loss-func}
We discuss four different choices of loss function $\rho(\cdot)$.
These loss functions are summarized in Table \ref{tab:loss-and-weight-functions}.
For each loss function, we also summarize  the  consistency factor $b$ and  the weight function  $u(t) = \mathrm{d} \rho(t)/\mathrm{d}t $  associated with the loss function $\rho$. 

\subsubsection{Gauss loss} corresponds to loss function of (circular complex) Gaussian distribution:
\begin{equation}\label{eq:loss-G}
\rho_{\mathrm{Gauss}}(t)=t
\end{equation}
in which case the consistency factor $b=1$ and the objective in \eqref{eq:Mobjective} becomes the Gaussian (negative) log-likelihood function $\mathrm{tr}\{ \Mat{\Sigma}^{-1} \Mat{S}_{\Mat{Y}} \}  - \log\det(\Mat{\Sigma}^{-1}) $ where 
\begin{align}
\Mat{S}_{\Mat{Y}} = \Mat{YY}^{\sf H} / L,
\label{eq:SY}    
\end{align}
is the SCM,  $(\cdot)^{\sf H}$ denotes Hermitian transpose,  and the minimizer of which is $\hat{\Mat{\Sigma}} = \Mat{S}_{\Mat{Y}}$.  In this case $b$ in \eqref{eq:consitency factor2} becomes $b=1$ as expected, since $\Mat{S}_{\Mat{Y}}$ is consistent to $\Mat{\Sigma}$ without any scaling correction. 
For Gauss loss $u_{\mathrm{Gauss}}(t)=1$. 

\subsubsection{Huber loss} given by \cite[Eq.~(4.29)]{Zoubir2018}
\begin{equation} \label{eq:loss-H} 
\rho_{\mathrm{Huber}}(t;c) =   \begin{cases} t &   \ \mbox{for} \ t \leqslant c^2, \\ 
c^2 \big( \log (t/c^2) + 1 \big)  & \ \mbox{for} \ t > c^2. \end{cases}
\end{equation} 
The threshold $c$ is a tuning parameter that affects the robustness and efficiency of the estimator.  Huber loss specializes the objective function \eqref{eq:Mobjective} to the negative log-likelihood of $\Mat{\Sigma}$ when the array data are heavy-tailed CES distributed with a density generator of the form $\mathrm{e}^{-\rho_{\mathrm{Huber}}(t;c)}$.  The squared threshold $c^2$ 
in \eqref{eq:loss-H} is mapped to the  $q$th quantile of $(1/2) \chi^2_{2N}$-distribution and we regard $q \in (0,1)$ as a loss parameter which is chosen by design, see Table \ref{tab:loss-and-weight-functions}.  
It is easy to verify that $b$ in \eqref{eq:consitency factor2}  for Huber loss function is  \cite[Sec.~4.4.2]{Zoubir2018},
\begin{align}\label{eq:consistency-H}
b_{\mathrm{Huber}} &= F_{\chi^2_{2(N+1)}}(2c^2) + c^2(1- F_{\chi^2_{2 N}}(2 c^2))/N, \\
   &= F_{\chi^2_{2(N+1)}}(2c^2) + c^2(1- q)/N, 
\end{align}
where $F_{\chi^2_{2N}}(x)$ denotes the cumulative distribution of the $\chi^2_{2N}$ distribution.  

For Huber loss \eqref{eq:loss-H} the weight function becomes 
\be \label{eq:huber_weight}
u_{\mathrm{Huber}} (t;c) = \begin{cases}  1,  
&  \ \mbox{for} \ t \leqslant  c^2 \\ c^2/t ,  & \ \mbox{for} \ t > c^2 \end{cases} . 
\ee 
Thus, an observation $\Vec{y}_{\ell}$ with squared Mahalanobis distance (MD) $\Vec{y}^{\sf H}_{\ell} \Mat{\Sigma}^{-1} \Vec{y}_{\ell}$ smaller than $c^2$ receives constant weight, while observations with a larger MD are heavily down-weighted.  

\subsubsection{MVT loss} which corresponds to the ML-loss for  (circular complex) multivariate $t$ (MVT) distribution with $\nu_{\mathrm{loss}}$ degrees of freedom, $\Vec{y}_i \sim \mathbb{C} t_{N,\nu}(\Vec{0},\Mat{\Sigma})$  \cite[Eq.(4.28)]{Zoubir2018}, 
\begin{equation}  \label{eq:loss-T}
\rho_{\mathrm{MVT}}(t;\nu_{\mathrm{loss}})=  \frac{\nu_{\mathrm{loss}} + 2N}{2} \log (\nu_{\mathrm{loss}} +2 t).
\end{equation} 
The $\nu_{\mathrm{loss}}$ parameter in \eqref{eq:loss-T} is viewed as a loss parameter which is chosen by design,   see Table \ref{tab:loss-and-weight-functions}.  The consistency factor $b_{_\mathrm{MVT}}$ 
for $\rho_{\mathrm{MVT}}(t;\nu_{\mathrm{loss}})$ is computable by numerical integration. 
 
 For MVT-loss  \eqref{eq:loss-T} the corresponding weight function is 
 \be \label{eq:t_weight}
u_{\mathrm{MVT}}(t; \nu_{\mathrm{loss}})=\frac{\nu_{\mathrm{loss}}+2N}{\nu_{\mathrm{loss}} + 2t},
\ee

\subsubsection{Tyler loss} given by \cite{Tyler1987jstor}, \cite[Sec. 4.4.3, Eq.  (4.30)]{Zoubir2018}
\begin{equation} \label{eq:loss-Tyler}
\rho_{\mathrm{Tyler}}(t)=  N \log (t).
\end{equation}
which is the limiting case of $\rho_{\mathrm{MVT}}(t;\nu_{\mathrm{loss}})$ for $\nu_{\mathrm{loss}}\to0$  and 
of $\rho_{\mathrm{Huber}}(t;c)$ for $c\to0$.  
This is the ML-loss of the Angular Central Gaussian (ACG) distribution \cite[Sec.  4.2.3]{Zoubir2018} which is not a CES distribution.
For Tyler loss \eqref{eq:loss-Tyler} the weight function becomes 
\begin{equation}
u_{\mathrm{Tyler}} (t) = N/t.
\end{equation}
In this case,  Tyler's M-estimator estimates the shape of the covariance matrix $\boldsymbol{\Sigma}$ only.  Namely, for Tyler loss and  $b=1$,  if $\hat{\boldsymbol{\Sigma}}$ is a minimizer of \eqref{eq:Mobjective}  then so is $b \hat{\boldsymbol{\Sigma}}$ for any $b > 0$. Thus,  the solution is unique only up to a scale. 
 
For Tyler loss $\psi(t)\equiv N$ $\forall t$,  indicating that the consistency factor for Tyler loss can not be found based on \eqref{eq:consitency factor}. 
In this case,  $\Mat{R}_{\Mat{Y}}$ is defined as weighted SCM that uses normalized Tyler weights \eqref{eq:consistency-factor-Tyler},  see Appendix \ref{sec:consistency-factor}.

\begin{table}
\begin{center}
\small
\begin{tabular}{|l|c|c|c|l|} \hline
 loss     &  loss      & weight         & &loss  \\ 
 name  & $\rho(t)$ & $u(t)$ & $\psi(t)$ &  parameter  \\ \hline
Gauss & \eqref{eq:loss-G} & 1 & $t$ & n/a \\
MVT &   \eqref{eq:loss-T} & \eqref{eq:t_weight} &
                                               $ \frac{\nu_{\mathrm{loss}}+2N}{2+\nu_{\mathrm{loss}}/t}$
                                                                                             & $\nu_{\mathrm{loss}}=2.1$  \\
Huber &  \eqref{eq:loss-H}  & \eqref{eq:huber_weight} & $t\,u_{\mathrm{Huber}}(t)$ & $q=0.9$  \\ 
Tyler &  \eqref{eq:loss-Tyler} & $N/t$ & $N$ & n/a \\ \hline
\end{tabular}
\\*[0.5ex]
\caption{Loss and weight functions used in DOA M-estimation with their loss parameter.
\label{tab:loss-and-weight-functions}}
\end{center}
\end{table}

\subsection{Source Power Estimation}
Similarly to Ref. \cite[Sec. III.D]{gerstoft2016mmv}, we  regard 
\eqref{eq:Mobjective} as a function of $\Vec{\gamma}$ and $\sigma^2$ and compute the first order derivative
\begin{align}
    \frac{\partial\mathcal{L}}{\partial \gamma_m} &= -\Vec{a}_m^{\sf H} \Mat{\Sigma}^{-1} \Vec{a}_m 
    + \frac{1}{Lb} \sum\limits_{\ell=1}^L  \| \Vec{a}_m^{\sf H} \Mat{\Sigma}^{-1} \Vec{y}_{\ell} \|_2^2 u(\Vec{y}_{\ell}^{\sf H} \Mat{\Sigma}^{-1} \Vec{y}_{\ell};\cdot) .
    \label{eq:dL-over-dgamma-m-new}  
 \end{align}
Equation \eqref{eq:dL-over-dgamma-m-new} is identical to  Ref. \cite[Eq.(21)]{gerstoft2016mmv} except for the weight function $u(\Vec{y}_{\ell}^{\sf H} \Mat{\Sigma}^{-1} \Vec{y}_{\ell};\cdot)$.  For the Gaussian array data model,  the weight function is the constant function $u_{\mathrm{Gauss}}(t)\equiv1$.

Setting \eqref{eq:dL-over-dgamma-m-new}  to zero gives
\begin{align}
    \Vec{a}_m^{\sf H} \Mat{\Sigma}^{-1} \Vec{a}_m 
    &= \Vec{a}_m^{\sf H} \Mat{\Sigma}^{-1} \Mat{R}_{\Vec{Y}} 
     \Mat{\Sigma}^{-1} \Vec{a}_m,  \; \forall m\in\{1,\ldots,M\}, 
     \label{eq:Jaffer-new}
\end{align}
where $\Mat{R}_{\Vec{Y}}$ is the weighted SCM
\begin{align}
\Mat{R}_{\Vec{Y}} &=\frac{1}{Lb}
\sum\limits_{\ell=1}^L u(\Vec{y}_{\ell}^{\sf H} \Mat{\Sigma}^{-1} \Vec{y}_{\ell};\cdot)
\Vec{y}_{\ell} \Vec{y}_{\ell}^{\sf H} \label{eq:RY} 
=\frac{1}{L}\Mat{Y}\Mat{D}\Mat{Y}^{\sf H}
\end{align}
with $\Mat{D} = \mathrm{diag}(u_1,\ldots,u_L)/b$ and $u_{\ell} = u(\Vec{y}_{\ell}^{\sf H} \Mat{\Sigma}^{-1} \Vec{y}_{\ell};\cdot)$.
Note that $\Mat{R}_{\Vec{Y}}$ can be understood as an adaptively weighted SCM \cite[Sec. 4.3]{Zoubir2018}.
$\Mat{R}_{\Vec{Y}}$ is Fisher consistent for the covariance matrix when $\Vec{y}_\ell$ are Gaussian distributed,
i.e. ,  $\mathsf{E}[\Mat{R}_{\Vec{Y}}]=\Mat{\Sigma}$ thanks to the consistency factor $b$ \cite[Sec. 4.4.1]{Zoubir2018}.
Although \eqref{eq:Jaffer-new} is rather concise,  it hides the parameters $\Vec{\gamma},  \sigma^2$ in $\Mat{\Sigma}$ and $\Mat{R}_{\Mat{Y}}$ and a closed-form solution for $\Vec{\gamma},  \sigma^2$ is not known to the authors.  Therefore,  we formulate a fixed-point equation for each $\gamma_m$ given $\sigma^2$ and $\gamma_p$, $p\ne m\in\{1,\ldots,M\}$ for solving \eqref{eq:Jaffer-new} by iterations.  Let $\gamma_m^{\mathrm{old}}$ be a previous approximation for $\gamma_m$ then (similarly to EM-based approaches \cite[Eqs.  (18)--(19)]{Wipf2007} and \cite{Nannuru2018}) the gradient-based update is computed with stepsize $\mu$,  cf. \cite[Ch.  8]{Sayed2008}
\begin{align}
     \gamma_m^{\mathrm{new}} &=  (1-\mu)   \gamma_m^{\mathrm{old}} + \mu  \gamma_m^{\mathrm{old}} G_m(\Vec{\gamma}^{\mathrm{old}}),  \quad\forall m, \label{eq:gamma-update-rule-with-stepsize}
\end{align}
where $0 < \mu \le 1$ and
\begin{align}
      G_m(\Vec{\gamma}) &=  \frac{\Vec{a}_m^{\sf H} \Mat{\Sigma}^{-1} \Mat{R}_{\Vec{Y}} \Mat{\Sigma}^{-1} \Vec{a}_m}{\Vec{a}_m^{\sf H}\Mat{\Sigma}^{-1} \Vec{a}_m}, \quad\forall m. \label{eq:iteration-gain-Gm} \\
 &= \frac{\frac{1}{L} \sum\limits_{\ell=1}^L\left| \Vec{a}_m^{\sf H} \Mat{\Sigma}^{-1}\Vec{y}_{\ell}\sqrt{u_{\ell}/b}\right|^2}{\Vec{a}_m^{\sf H}\Mat{\Sigma}^{-1} \Vec{a}_m}.  
\end{align}
The definition of $G_m(\Vec{\gamma})$ in \eqref{eq:iteration-gain-Gm} ensures that \eqref{eq:Jaffer-new} is fulfilled at the fixed-point of \eqref{eq:gamma-update-rule-with-stepsize}.
In \eqref{eq:iteration-gain-Gm},  both $\Mat{\Sigma}^{-1}$ and $\Mat{R}_{\Mat{Y}}$ depend on $\Vec{\gamma}$.
The simulations results in Sec.  \ref{sec:results} are obtained for stepsize $\mu=1$.
The stability of iterations based on \eqref{eq:gamma-update-rule-with-stepsize} is discussed in Appendix \ref{sec:stability} and numerical observations of the associated run-time are reported in Sec. \ref{sec:convergence}.

The active set $\mathcal M$ is selected as either the $K$ largest entries of  $\Vec{\gamma}$ or the entries with $\gamma_m$ exceeding a threshold.

\subsection{Noise Variance Estimation}

The original SBL algorithm exploits Jaffer's necessary condition \cite[Eq.~(6)]{Jaffer1988} which leads to the noise subspace based estimate \cite[Eq.~(15)]{Liu2012}, \cite[Sec. III.E]{gerstoft2016mmv},
\begin{align}
\hat{\sigma}^2_S &= 
 \frac{\trace{(\Mat{I}_N-\Mat{A}^{}_{\mathcal{M}} \Mat{A}^+_{\mathcal{M}} )\Mat{S}_{\Mat{Y}}}}{N-K } ,
\label{eq:noise-estimate}
\end{align}
where $(\cdot)^+$ denotes the Moore-Penrose pseudo inverse.
This noise variance estimate works well with DOA estimation~\cite{gerstoft2016mmv,Nannuru2018,Park2021} without outliers in the array data. 
The derivation of the robust noise variance estimate starts from the following robust formulation of Jaffer’s condition, namely
\begin{align}
        \Mat{A}_{\mathcal{M}}^{\sf H} (\Mat{R}_{\Mat{Y}}-\Mat{\Sigma}) \Mat{A}_{\mathcal{M}} = \Mat{0}.
\end{align}
This results in the robust estimate
\begin{align}
\hat{\sigma}^2_R &= 
 \frac{\trace{(\Mat{I}_N-\Mat{A}^{}_{\mathcal{M}} \Mat{A}^+_{\mathcal{M}})\Mat{R}_{\Mat{Y}}}}{N-K } ,
\label{eq:noise-estimateR}
\end{align}
for the full derivation see \cite{Mecklenbraeuker-Gerstoft-Ollila-WSA2021}. For Gauss loss function, $\Mat{R}_{\Mat{Y}}=\Mat{S}_{\Mat{Y}}$, and the expressions \eqref{eq:noise-estimate}, \eqref{eq:noise-estimateR} are identical.

To stabilize the noise variance M-estimate \eqref{eq:noise-estimateR} for 
non-Gauss loss,  we define lower and upper bounds for $\hat{\sigma}^2$
and enforce
$\sigma^2_{\text{floor}}\le \hat{\sigma}^2\le \sigma^2_{\text{ceil}}$ by
 \begin{equation}
 \hat{\sigma}^2=\max(\min(\hat{\sigma}^2_R,\sigma^2_{\text{ceil}}), \sigma^2_{\text{floor}})
\label{eq:sigma-stabilization}
 \end{equation}
The original SBL algorithm \cite{gerstoft2016mmv,SBL4-github} does not use/need this stabilization because for Gauss loss,  the weighted SCM \eqref{eq:RY} equals \eqref{eq:SY} which does not depend on prior knowledge of $\Mat{\Sigma}$.
We have chosen $\sigma^2_{\text{floor}}=10^{-6}{\rm tr}[\Mat{S}_{\Mat{Y}}]/N$ and $\sigma^2_{\text{ceil}}={\rm tr}[\Mat{S}_{\Mat{Y}}]/N$ for the numerical simulations.

As discussed in Appendix \ref{sec:consistency-factor}, Tyler's M-estimator is unique only up to a scale which affects the noise variance estimate $\hat{\sigma}^2_R$.  For this reason,  we normalize $\Mat{R}_{\Mat{Y}}$ to trace 1  to remove this ambiguity if Tyler loss is used.

\begin{table}
\begin{algorithmic}[1]
  \small
  \STATE {\bf input} $\Mat{Y}\in\mathbb{C}^{N\times L}$ \hspace{15ex}  array data to be analyzed
  \STATE select the weight function $u(\cdot;\cdot)$ and loss parameter
  \STATE constant $ \Mat{A}\in\mathbb{C}^{N\times M}$ \hspace{12ex}  dictionary matrix
  \STATE \hspace{8ex} $K\in\mathbb{N}$, with $K<N$,  \hspace{2.5ex}  number of sources
  \STATE \hspace{8ex} $\delta\in\mathbb{R}^+$ \hspace{16ex}  small positive constant
  \STATE \hspace{8ex} $\text{SNR}_{\text{max}}\in\mathbb{R}^+$ \hspace{10ex} upper SNR limit in data
  \STATE \hspace{8ex} $\gamma_{\text{range}}\in[0,1]$,  \hspace{0.5ex}  dynamic range for DOA grid pruning
  \STATE \hspace{8ex} $\mu\in[0,1]$, \hspace{14ex} iteration stepsize
  \STATE \hspace{8ex} $j_{\text{max}}\in\mathbb{N}$  \hspace{15ex} iteration count limit 
  \STATE \hspace{8ex} $z\in\mathbb{N}$ with $z<j_{\text{max}}$ \hspace{4.2ex} convergence criterion 
  \STATE  set $\Mat{S}_{\Mat{Y}} = \Mat{Y} \Mat{Y}^{\sf H}/L$
  \STATE  set $\sigma^2_{\text{ceil}} = \trace{\Mat{S}_{\Mat{Y}}}$ \hspace{15ex} upper limit on $\hat{\sigma}^2$
  \STATE  \hspace{2.5ex} $\sigma^2_{\text{floor}} = \sigma^2_{\text{ceil}}/\text{SNR}_{\text{max}}$ \hspace{9ex} lower limit on $\hat{\sigma}^2$
  \STATE initialize $\hat{\sigma}^2$, $\Vec{\gamma}^{\text{new}}$ using \eqref{eq:init-step1}--\eqref{eq:init-step4},  
                                $j=0$ 
  \REPEAT
  \STATE $j = j + 1$ \hspace{20ex} increment iteration counter
  \STATE $\Vec{\gamma}^{\text{old}}=\Vec{\gamma}^{\text{new}}$
  \STATE $\gamma_{\text{floor}} = \gamma_{\text{range}} \max(\Vec{\gamma}^{\text{new}})$,  \hspace{6ex} source dynamic range
  \STATE $\mathcal{P}=\{p \in \mathbb{N} \; \vert \;   \Vec{\gamma}^{\text{new}}_p \ge \gamma_{\text{floor}}\}$ \hspace{3ex} pruned DOA grid \hfill \eqref{eq:mathcal-P}
  \STATE $\Mat{\Gamma}_{\mathcal{P}}=\text{diag}(\Vec{\gamma}^{\text{new}}_{\mathcal{P}})$ \hspace{13ex} pruned source powers
  \STATE $\Mat{A}_{\mathcal{P}}=[\Vec{a}_{p_1}, \dots, \Vec{a}_{p_P}]$ for all $p_i\in\mathcal{P}$,  pruned dictionary
  \STATE  $\Mat{\Sigma}_{\mathcal{P}}= \Mat{A}_{\mathcal{P}} \Mat{\Gamma}_{\mathcal{P}} \Mat{A}^{\sf H}_{\mathcal{P}} + \hat{\sigma}^2 \Mat{I}_N $ \hfill \eqref{eq:Sigma-model_new}
 \STATE {\bf if} using $u_{\mathrm{Tyler}}(\cdot)$ {\bf then} $b=\frac{1}{L}\sum_{\ell=1}^L u_{\mathrm{Tyler}}(\Vec{y}_\ell^{\sf H} \Mat{\Sigma}^{-1} \Vec{y}_\ell)$ \hfill \eqref{eq:consistency-factor-Tyler}
  \STATE $\Mat{R}_{\Vec{Y}} =\frac{1}{Lb}
\sum\limits_{\ell=1}^L u(\Vec{y}_{\ell}^{\sf H} \Mat{\Sigma}_{\mathcal{P}}^{-1} \Vec{y}_{\ell};\cdot)
\Vec{y}_{\ell} \Vec{y}_{\ell}^{\sf H}$ \hfill \eqref{eq:RY}
 \STATE $\gamma_p^{\mathrm{new}} =  (1-\mu)   \gamma_p^{\mathrm{old}} + \mu  \gamma_p^{\mathrm{old}} G_p(\Vec{\gamma}^{\mathrm{old}}),    $ for all $p\in\mathcal{P}$ \hfill \eqref{eq:gamma-update-rule-with-stepsize}
  \STATE $\mathcal{M}=\{m \in \mathbb{N} \; \vert \; K\; \text{largest peaks in } \Vec{\gamma}^{\text{new}} \}$ active set
  \STATE $\Mat{A}_{\mathcal{M}}=[\Vec{a}_{m_1}, \dots, \Vec{a}_{m_K}]$
  \STATE $\hat{\sigma}^2_R = \frac{\trace{(\Mat{I}_N- \Mat{A}^{}_{\mathcal{M}}\Mat{A}_{\mathcal{M}}^+)\Mat{R}_{\Mat{Y}}}}{N-K } $ \hfill \eqref{eq:noise-estimateR}
 \STATE $\hat{\sigma}^2=\max(\min(\hat{\sigma}^2_R,\sigma^2_{\text{ceil}}), \sigma^2_{\text{floor}})$  \hfill \eqref{eq:sigma-stabilization}
    \UNTIL{($\mathcal{M}$ has not changed during last $z$ iterations) or $j > j_{\text{max}}$ }
  \STATE {\bf output} $\mathcal{M}$ and $j$ 
\end{algorithmic}
\caption{Robust and Sparse M-Estimation of DOA.  This table documents the algorithm formulated
for single-frequency array data.
The Matlab function \textcolor{black}{\tt SBL\_v5p11.m} \cite{RobustSBL-github} implements 
Robust and Sparse M-Estimation of DOA from multi-frequency array data. 
}
\label{table:algorithm}
\end{table}

\subsection{Algorithm}
\label{sec:algorithm}
The proposed DOA M-estimation algorithm using SBL is displayed in Table \ref{table:algorithm} with the following remarks: 

\subsubsection{DOA grid pruning}
\label{sec:grid-pruning}
To reduce numerical complexity in the iterations,  we introduce the pruned DOA grid $\mathcal{P}$ by not wasting  computational resources on those DOAs which are associated with source power estimates below a chosen threshold value $\gamma_{\text{floor}}$, i.e. ,  we introduce a thresholding operation on the $\Vec{\gamma}^{\text{new}}$ vector.  
The pruned DOA grid is formally defined as an index set, 
\begin{align}
   \mathcal{P} = \{p \in \{1,\ldots,M \}\; \vert \;   \Vec{\gamma}^{\text{new}}_p \ge \gamma_{\text{floor}}\} =
                            \{ p_1, \ldots,  p_P \}.
  \label{eq:mathcal-P}
\end{align}
where  $\gamma_{\text{floor}}=\gamma_{\text{range}} \max{   \Vec{\gamma}^{\text{new}}}$ and we have chosen $\gamma_{\text{range}} =10^{-3}$.

\subsubsection{Initialization} 
In our algorithm we need to give initial values of source signal powers $\Vec{\gamma}$ and the noise variance $\sigma^2$. 
The initial estimates are computed via following steps: 
\begin{enumerate}[a)]
\item Compute $\Mat{S}_{\Mat{Y}}$ and CBF output powers
\begin{align}
\gamma_m^{\text{init}} &=  \frac{ \Vec{a}^{\sf H}_m \mathbf{S}_Y \Vec{a}_m }{ \| \Vec{a}_m \|^4} , \quad \forall m = 1,\ldots,M.
\label{eq:init-step1}
\end{align}
\item Compute the initial active set by identifying $K$ largest peaks in the CBF output powers, 
\begin{align}
\mathcal{M} &= \{m \in \mathbb{N} \; \vert \; K\; \text{largest peaks in } \Vec{\gamma}^{\text{init}} \}
\label{eq:init-step2}
\end{align}
\item Compute the initial noise variance 
\begin{align}
\hat{\sigma}^2 &= \hat{\sigma}^2_S =  \frac{\trace{(\Mat{I}_N-\Mat{A}^{}_{\mathcal{M}} \Mat{A}^+_{\mathcal{M}}) \Mat{S}_{\Mat{Y}}}}{N-K }
\label{eq:init-step3}
\end{align}
\item Compute initial estimates of source powers: 
\begin{align}
\gamma_m^{\text{new}} &=  \max(\delta,  (\gamma_m^{\text{init}}  - \hat{\sigma}^2 )),
\quad \mbox{for $m=1,\ldots,M$}. 
\label{eq:init-step4}
 \end{align}
where $\delta>0$ is a small number, guaranteeing that all initial $\gamma_m^{\text{new}}$ are positive. 
 \end{enumerate} 

\subsubsection{Convergence Criterion} 
\label{sec:convergece-criterion}

The DOA Estimates returned by the iterative algorithm in Table \ref{table:algorithm} are obtained from the active set $\mathcal{M}$. Therefore,  the active set is monitored for changes in its elements to determine whether the algorithm has converged.  If $\mathcal{M}$ has not changed during the last $z\in\mathbb{N}$ iterations,  then the repeat-until loop (lines 14--29 in Table \ref{table:algorithm}) terminates.  Here $z$ is a tuning parameter which allows to trade off computation time against DOA estimate accuracy.  To ensure that the iterations always terminate,  the maximum iteration count is defined as $j_{\mathrm{max}}$ with 
$z<j_{\mathrm{max}}$.

\section{Simulation Results}
\label{sec:results}

The proposed DOA M-estimation algorithm using SBL is displayed in Table \ref{table:algorithm}.
Numerical simulations are carried out for evaluating the root mean squared error (RMSE) of DOA  versus array signal to noise ratio (ASNR) based on synthetic array data  $\Mat{Y}$.
Synthetic array data are generated for three scenarios with $K=1,\ldots,3$ incoming plane waves and corresponding DOAs as listed in Table \ref{tab:source-scenarios}.  
The source amplitudes $\Vec{x}_{\ell}$ in \eqref{eq:scale-mixture} 
are complex circularly symmetric zero-mean Gaussian.
The wavefield is modeled according to the scale mixture \eqref{eq:scale-mixture} and it is observed by a uniform linear array with $N=20$ elements at half-wavelength spacing. 
The dictionary $\Mat{A}$ consists of $M=18001$ replica vectors for the high-resolution DOA grid
$\theta_m = -90^{\circ} + (m-1)\delta,~\forall m=1,\ldots,M $ where $\delta$ is the dictionary's angular grid resolution, $\delta=180^\circ/(M-1)=0.01^{\circ}$.

The convergence criterion parameter $z=10$ is chosen for all numerical simulations and the maximum number of iterations was set to $j_{\mathrm{max}}=1200$, but this maximum was never reached.
The RMSE  of the DOA estimates  over $N_{\text{run}} = 250$ simulation runs with random array data realizations is used for evaluating the performance of the algorithm,
\begin{align}
\text{RMSE} &= \sqrt{\sum_{r=1}^{N_{\rm run}} \sum_{k=1}^{K} \frac{
[\min(|\hat{\theta}^r_k - \theta^r_k|,e_{\mathrm{max}})]^2}{K \,N_{\rm run}}} \,,
\label{eq:RMSE}
\end{align}
where $\theta^r_k$ is the true DOA of the $k$ source and $\hat{\theta}^r_k$ is the corresponding estimated DOA in the $r$th  run
when $K$ sources are present in the scenario. 
This RMSE definition is a specialization of the optimal subpattern assignment (OSPA) when $K$ is known,  cf.  \cite{Meyer2017}.
We use $e_{\mathrm{max}}=10^\circ$ in \eqref{eq:RMSE}. 
Thus,  maximum RMSE is $10^\circ$.

\subsection{Data generation}
The scaling $\sqrt{\tau}_{\ell}$ and the noise $\Vec{n}_{\ell}$ in \eqref{eq:scale-mixture} are generated according to three array data models which are summarized in Table \ref{tab:data-models} and explained below:

\paragraph{Gaussian array data} In this model $\tau_{\ell}=1$ for all $\ell$ in \eqref{eq:scale-mixture} and $\Vec{y}_{\ell}= \Vec{v}_{\ell}\sim\mathbb{C}\mathcal{N}(\Vec{0},\Mat{\Sigma})$, where $\Mat{\Sigma}$ is defined in \eqref{eq:Sigma-model_new}.

\paragraph{MVT array data} We first draw $\Vec{v}_{\ell}\sim\mathbb{C}\mathcal{N}(\Vec{0},\Mat{\Sigma})$ and $s_{\ell}\sim\chi^2_{\nu_{\mathrm{data}}}$ independently, where $\Mat{\Sigma}$ is defined in \eqref{eq:Sigma-model_new}.  We set
$\tau_\ell = \nu_{\mathrm{data}}/s_{\ell}$ and the array data is modelled by the scale mixture
$\Vec{y}_{\ell} = \sqrt{\tau_\ell}\Vec{v}_{\ell} \sim \mathbb{C}t_{\nu_{\mathrm{data}}}$-distributed, cf. \cite{Ollila2003pimrc} and \cite[Sec.~4.2.2]{Zoubir2018}.

\paragraph{$\epsilon$-contaminated array data} This heavy-tailed array data model is not covered by \eqref{eq:scale-mixture} with the assumptions in Sec. \ref{sec:model}.
Instead,  the noise $\Vec n$ is drawn with probability  $(1-\epsilon)$ from a $\mathbb{C}\mathcal{N}(\Vec{0}, \sigma^2_1\Mat{I})$ and with probability $\epsilon$ from a $\mathbb{C}\mathcal{N}(\Vec{0}, \lambda^2\sigma^2_1\Mat{I})$, where $\lambda$ is the outlier strength.
Thus, $\Vec{y}_{\ell}$ is drawn from
$\mathbb{C}\mathcal{N}(\Vec{0}, \Mat{A}\Mat{\Gamma}\Mat{A}^{\sf H}+\sigma^2_1\Mat{I}_N)$, using \eqref{eq:Sigma-model_new} with probability $(1-\epsilon)$ and with outlier probability $\epsilon$ from $\mathbb{C}\mathcal{N}(\Vec{0}, \Mat{A}\Mat{\Gamma}\Mat{A}^{\sf H}+(\lambda\sigma_1)^2\Mat{I}_N)$.
The resulting noise covariance matrix is $\sigma^2\Mat{I}_N$ similar to the other models, but with 
\begin{equation}
\sigma^2 = (1-\epsilon+\epsilon\lambda^2)\sigma_1^2.
\label{eq:epscont-variance}
\end{equation}
The limiting distribution of $\epsilon$-contaminated noise for $\epsilon \to 0$ and any constant $\lambda>0$ is Gaussian.

\begin{figure*}[t]
\parbox[t]{0.333\textwidth}{%
{\footnotesize a) \textcolor{white}{\tt p08fix\_Gmode\_s1MC250SNRn15.eps}}\\
\includegraphics[width=0.325\textwidth]{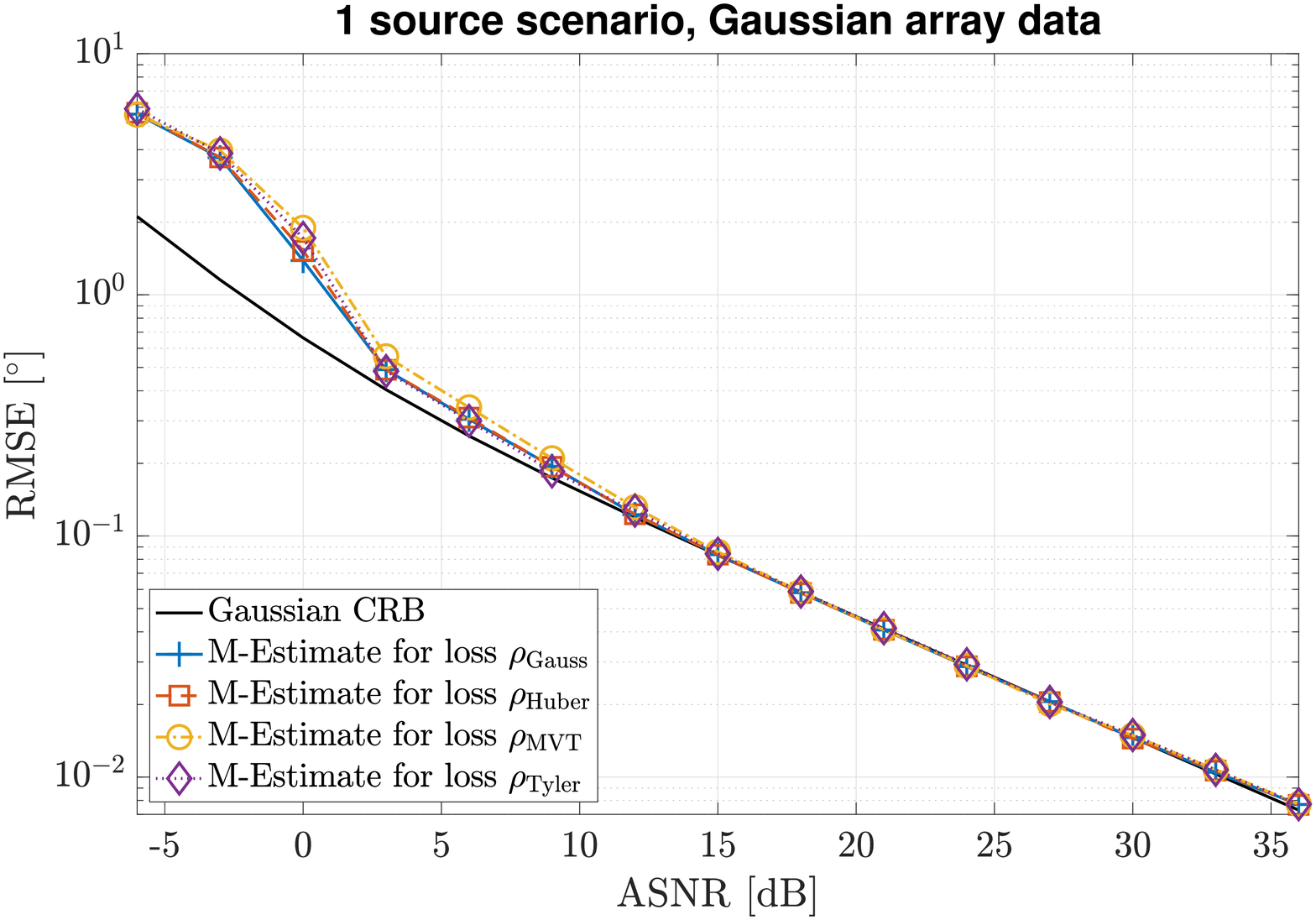}\\
{\footnotesize b) \textcolor{white}{\tt p08\_Cmode\_s1MC250SNRn15.eps}}\\
\includegraphics[width=0.325\textwidth]{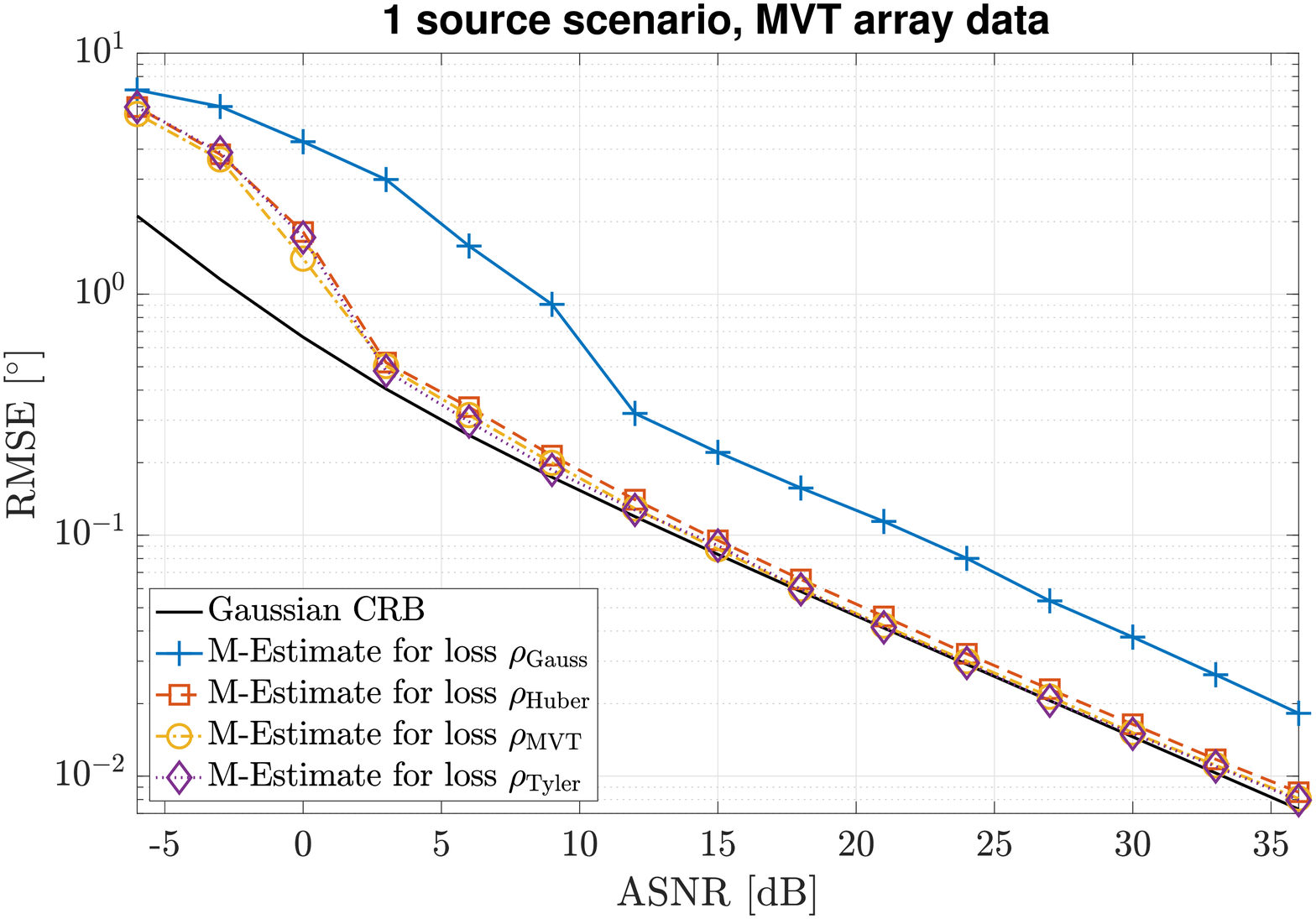}\\
{\footnotesize c) \textcolor{white}{\tt p08\_emode\_s1MC250SNRn15.eps}}\\
\includegraphics[width=0.325\textwidth]{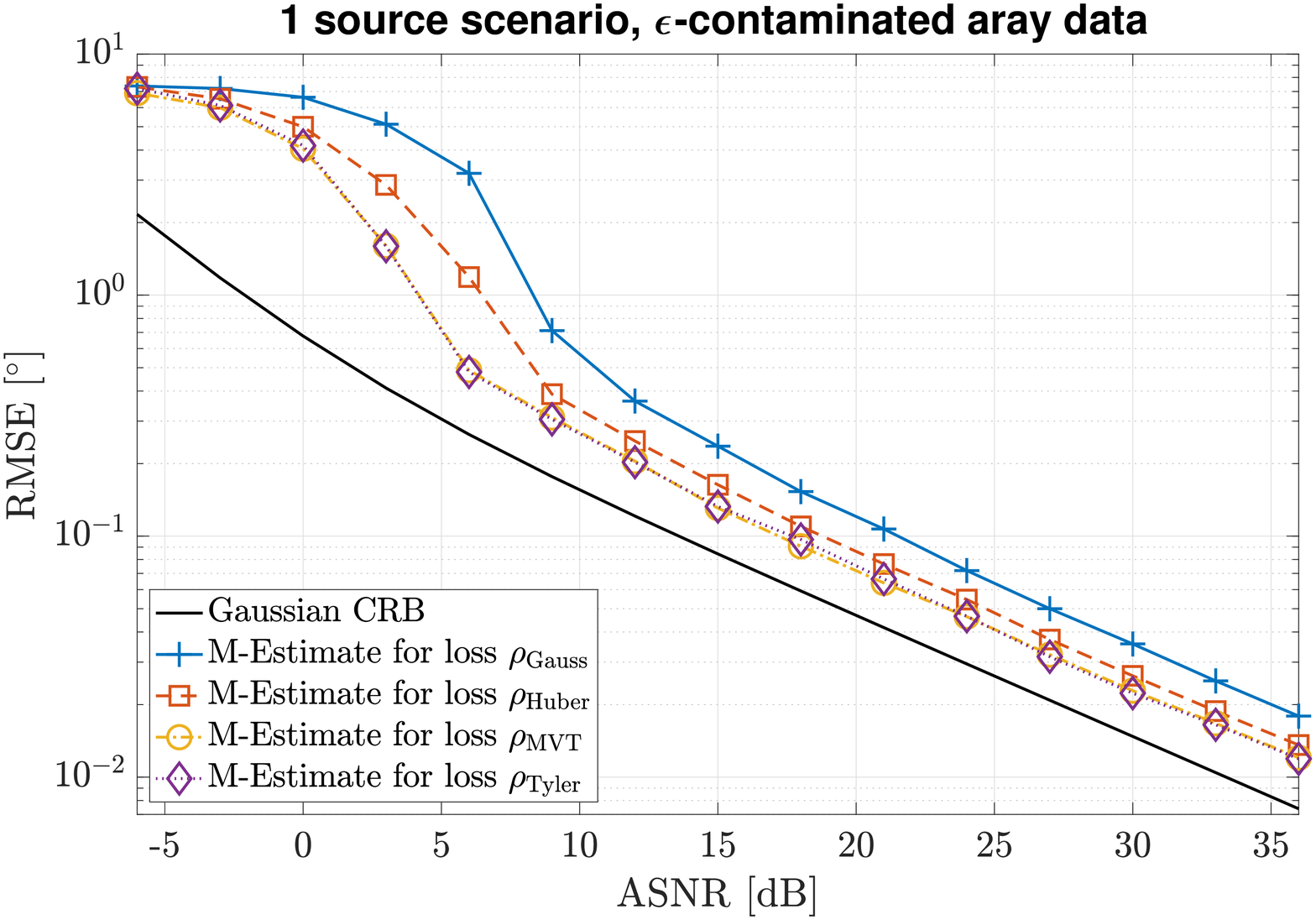}}
\parbox[t]{0.333\textwidth}{%
{\footnotesize d) \textcolor{white}{\tt p08\_Gmode\_s2MC250SNRn15.eps}}\\
\includegraphics[width=0.325\textwidth]{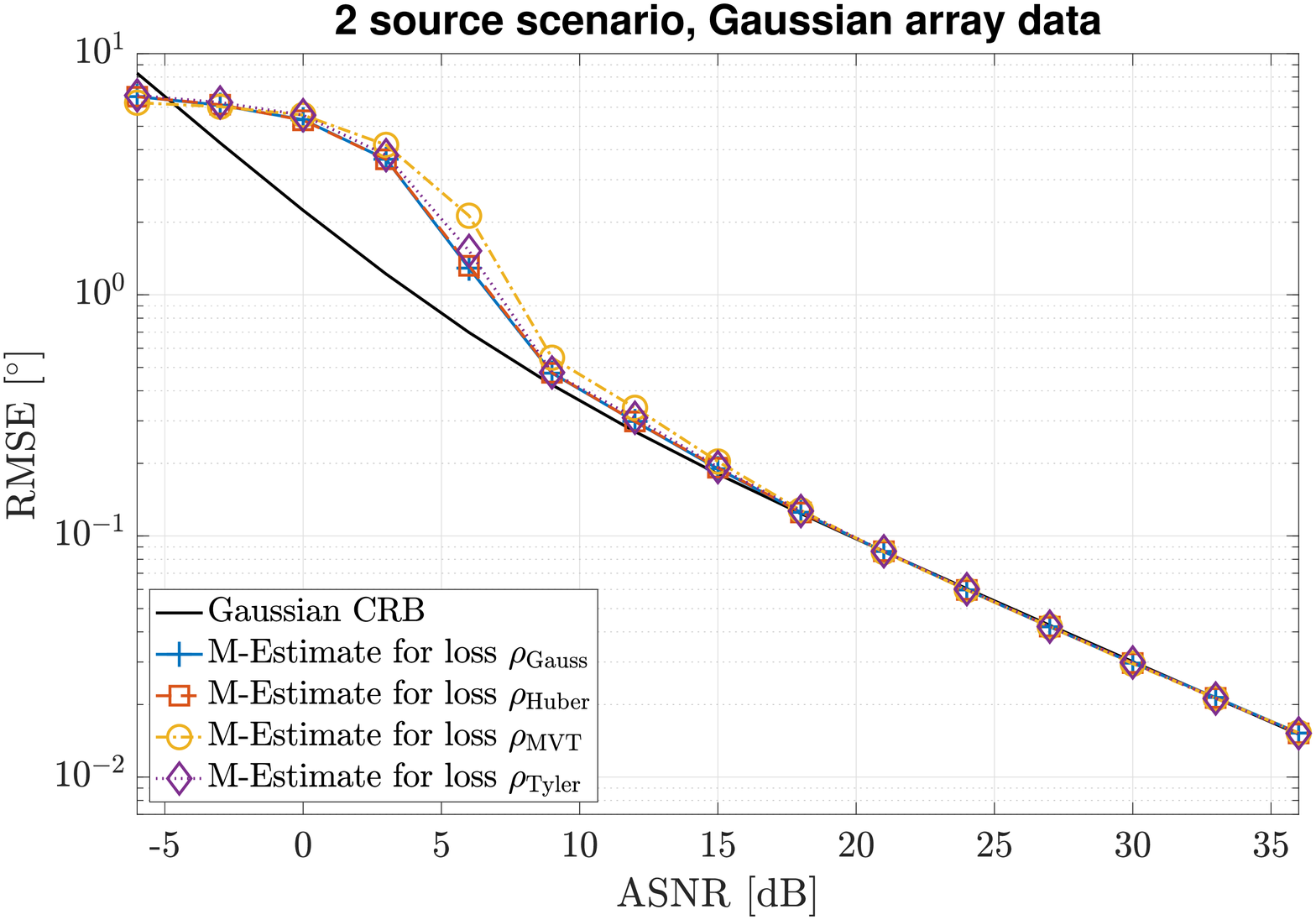}\\
{\footnotesize e) \textcolor{white}{\tt p08\_Cmode\_s2MC250SNRn15.eps}}\\
\includegraphics[width=0.325\textwidth]{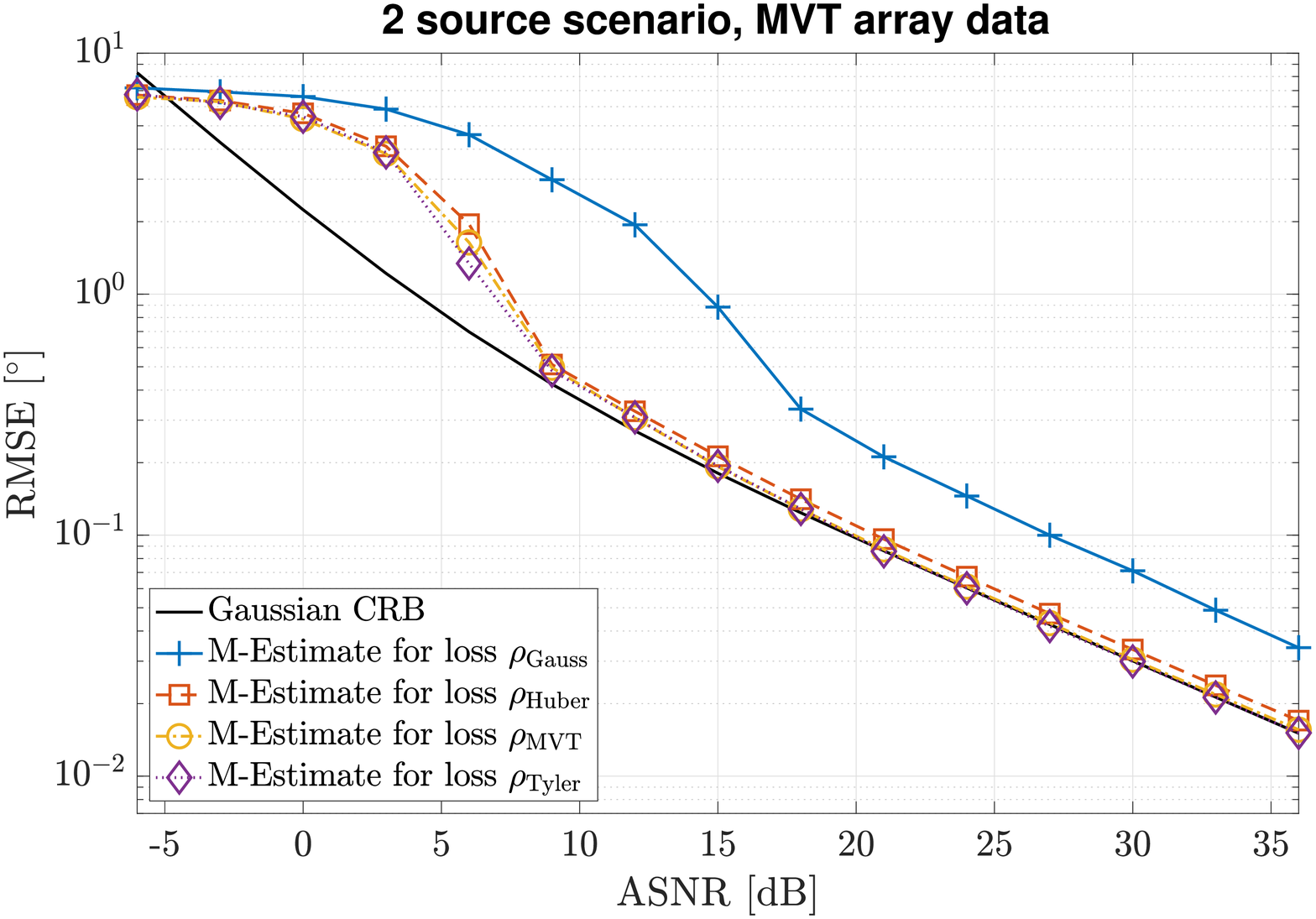}\\
{\footnotesize f) \textcolor{white}{\tt p08fix\_emode\_s2MC250SNRn15.eps}}\\
\includegraphics[width=0.325\textwidth]{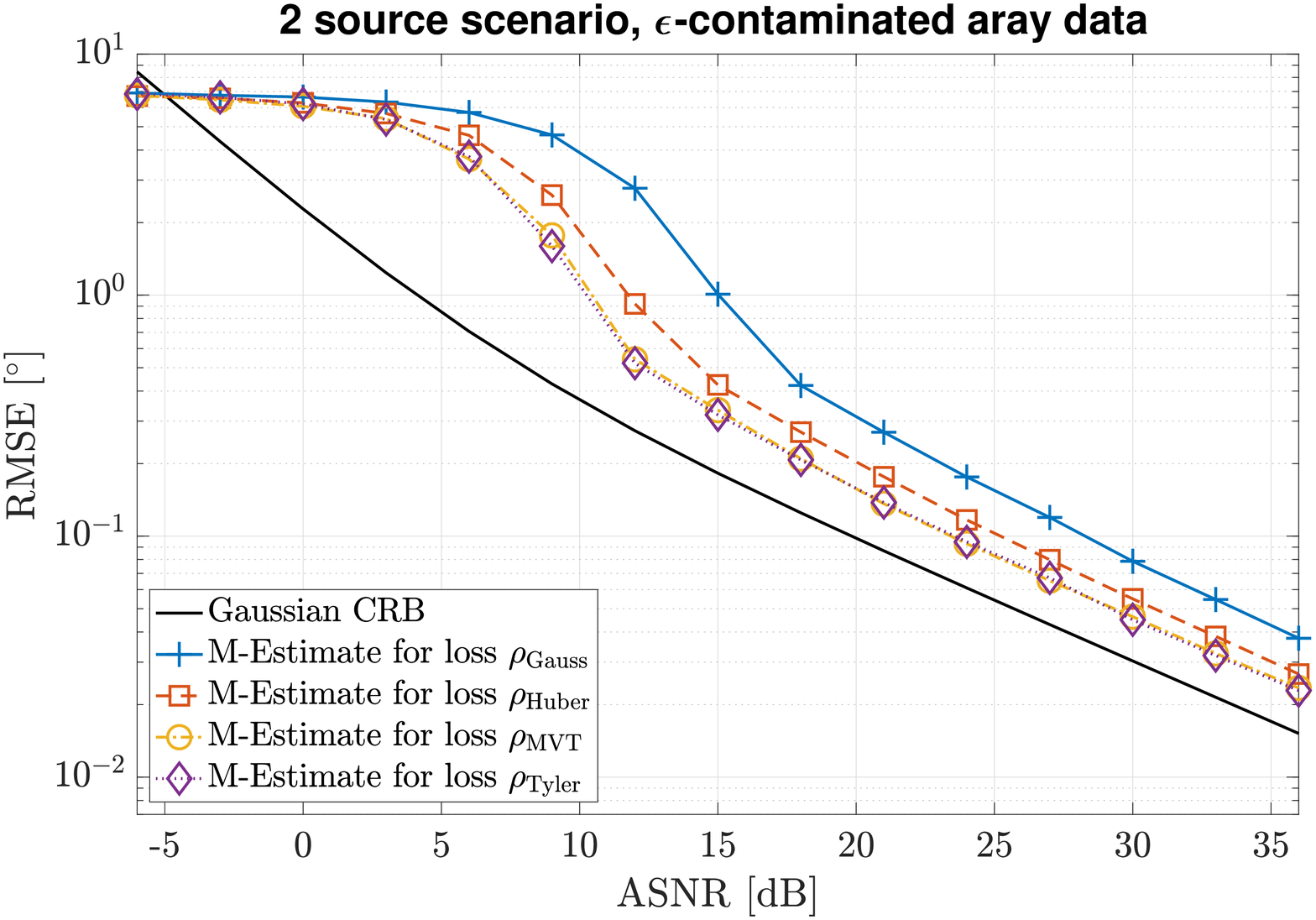}}
\parbox[t]{0.333\textwidth}{%
{\footnotesize g) \textcolor{white}{\tt p08fix\_Gmode\_s3MC250SNRn15.eps}}\\
\includegraphics[width=0.325\textwidth]{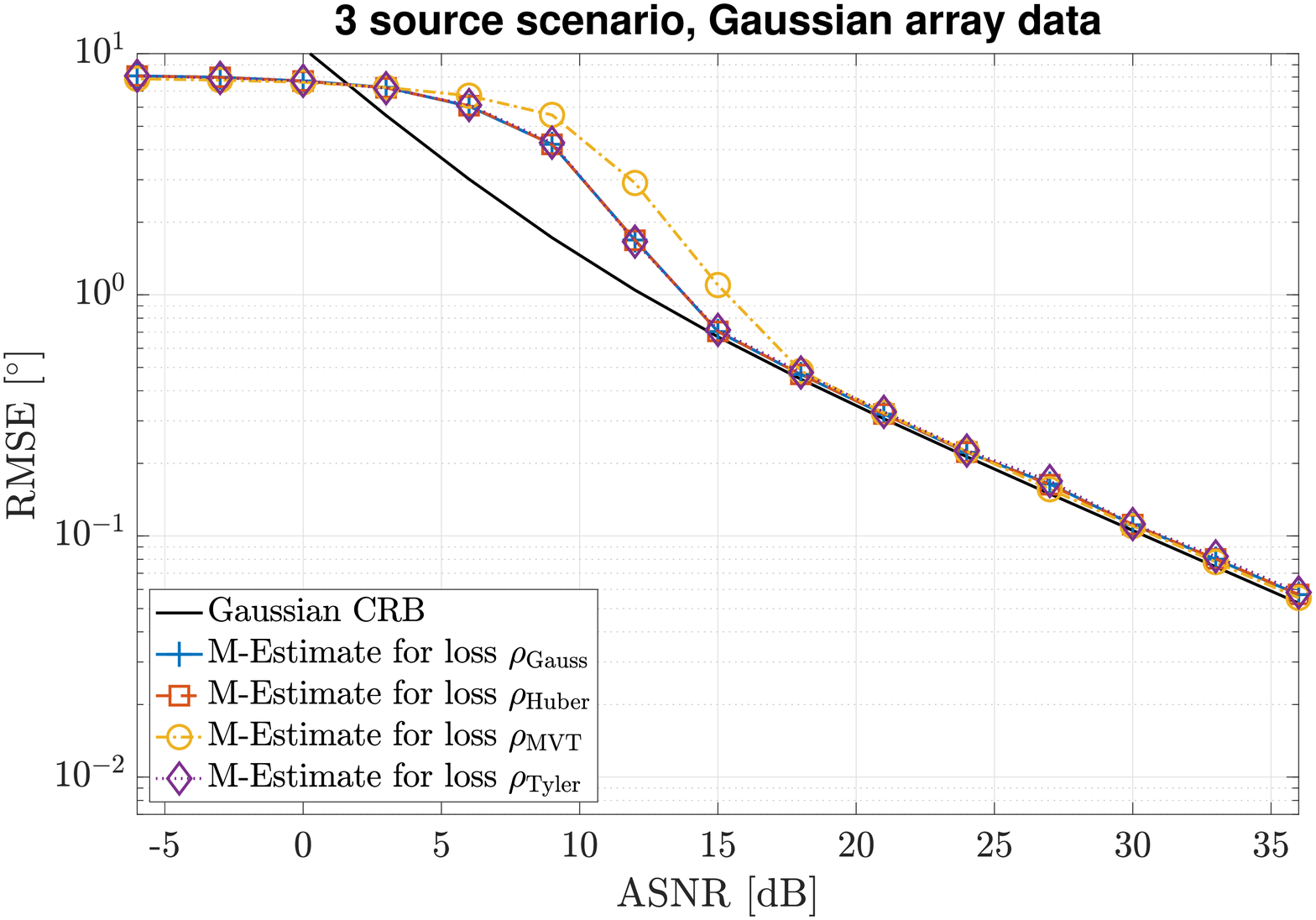}\\
{\footnotesize h) \textcolor{white}{\tt p08\_Cmode\_s3MC250SNRn15.eps}}\\
\includegraphics[width=0.325\textwidth]{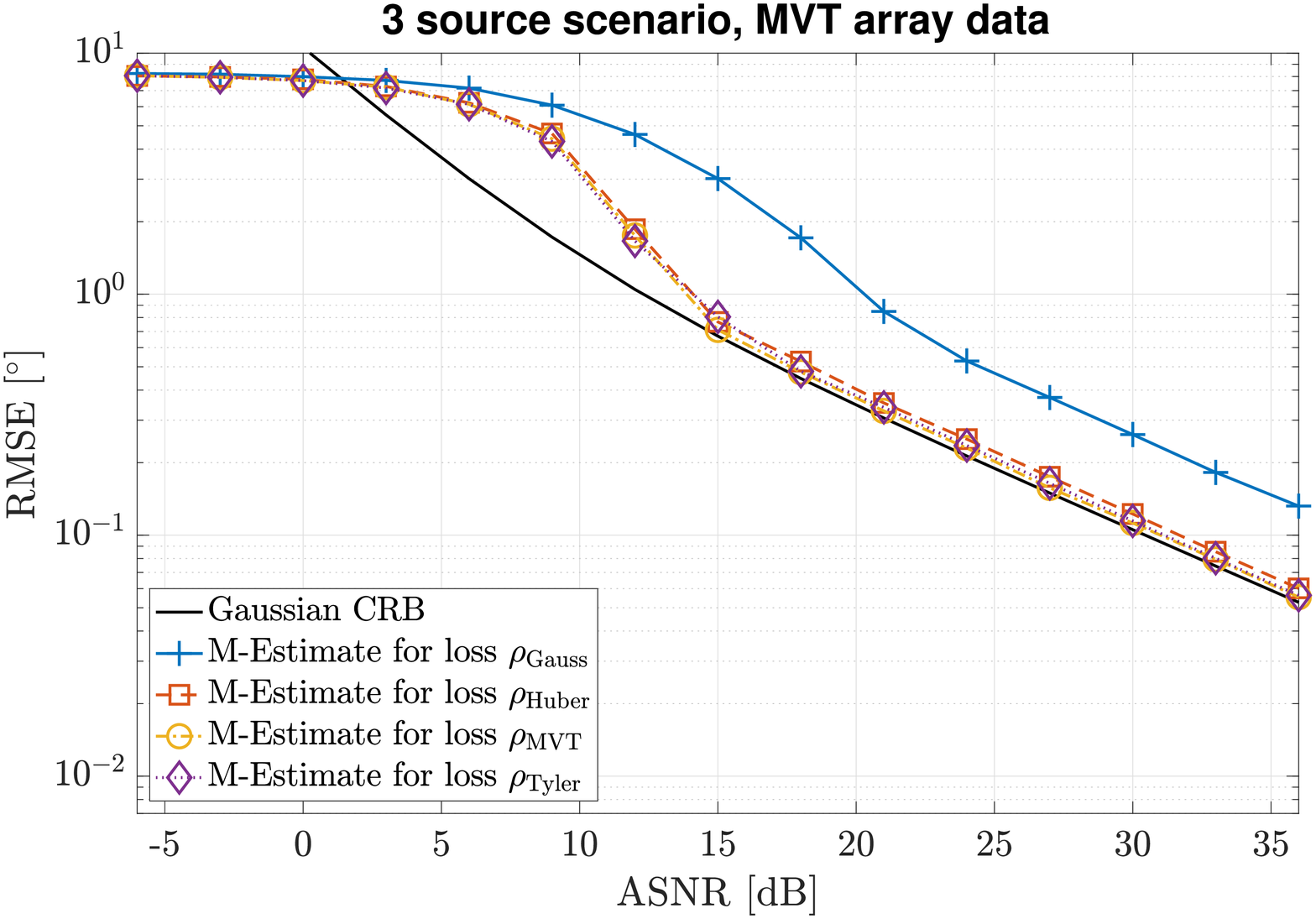}\\
{\footnotesize i) \textcolor{white}{\tt p08\_emode\_s3MC250SNRn15.eps}}\\
\includegraphics[width=0.325\textwidth]{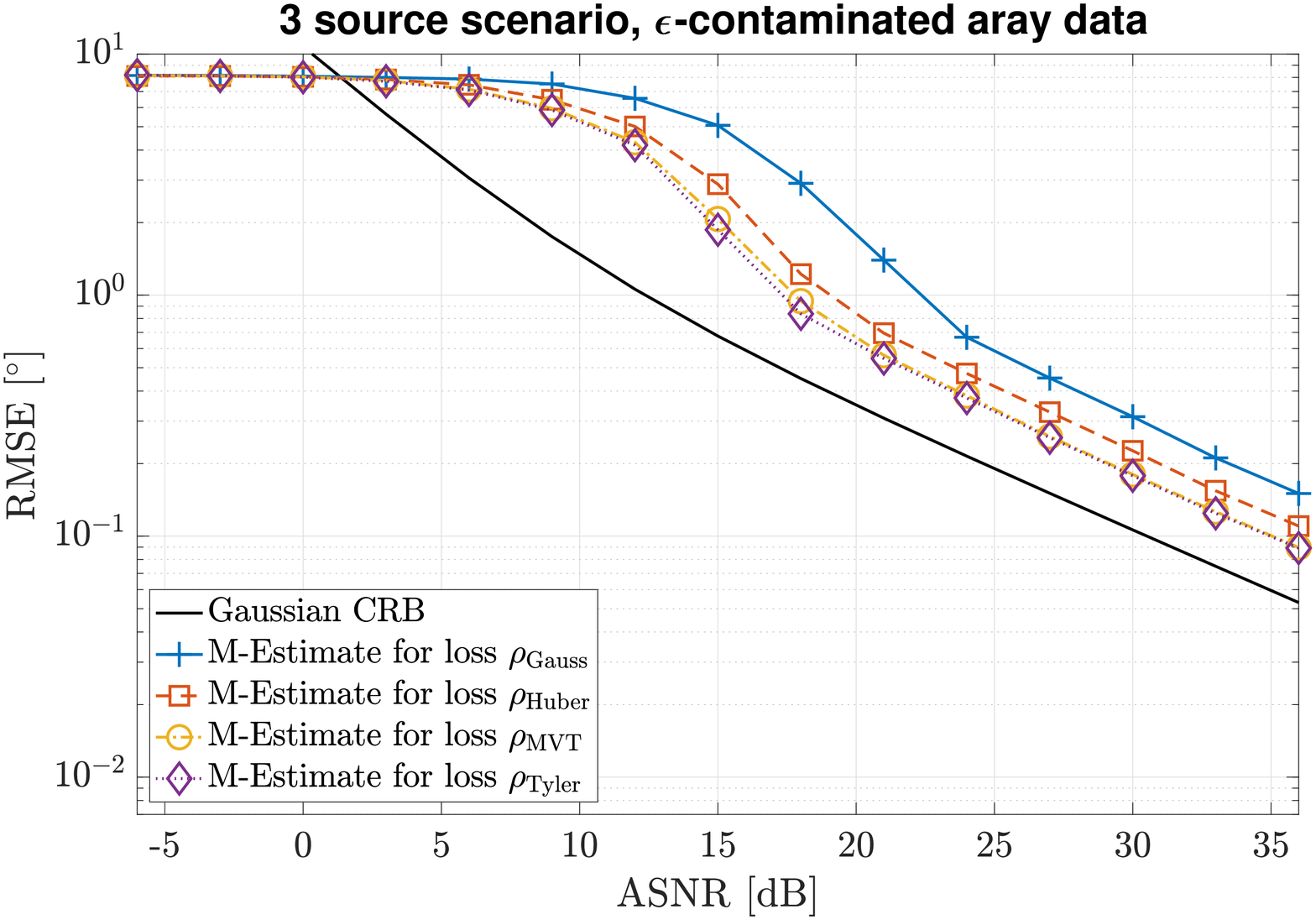}}
     \caption{RMSE of DOA estimators vs.~ASNR.  Left column: for single source at DOA $-10^{\circ}$.  Center column: for two sources at DOAs $-10^{\circ}$ and $10^{\circ}$. 
Right column: for three sources at DOAs $-3^{\circ}$,  $2^{\circ}$ and $75^{\circ}$. 
All: Simulation for uniform line array, $N=20$ sensors, $L=25$ array snapshots, and dictionary size $M=18001$ corresponding to DOA resolution $0.01^{\circ}$, averaged over $250$ realizations.  Array data models: Top row: Gaussian,  middle row: MVT ($\nu_{\mathrm{data}}=2.1$),  bottom row: $\epsilon$-contaminated   ($\epsilon=0.05,\lambda=10$). The CRB is for Gaussian \eqref{eq:Gaussian-CRB} and MVT array data models.,  see Appendix \ref{sec:CRB}. }
\label{fig:s1s2s3MC250SNRn15}
\end{figure*}

Additionally, the Cram\'er-Rao Bound (CRB) for DOA estimation from Gaussian and MVT array data are shown.  The corresponding expressions are given in Appendix \ref{sec:CRB} for completeness.
The Gaussian CRB shown in Figs.  \ref{fig:s1s2s3MC250SNRn15}(a,d,g) is evaluated according to \eqref{eq:Gaussian-CRB}.
The bound for MVT array data, $C_{\mathrm{CR,MVT}}(\Vec{\theta})$,  is just slightly higher than the bound for Gaussian array data,  $C_{\mathrm{CR,Gauss}}(\Vec{\theta})$.
For the three source scenario,  $N=20$,  $L=25$,  Gaussian and MVT array data model,  the gap between the bounds $C_{\mathrm{CR,Gauss}}(\Vec{\theta})$ and $C_{\mathrm{CR,MVT}}(\Vec{\theta})$ is smaller than $3\%$ in terms of RMSE in the shown ASNR range.  
The gap is not observable in the RMSE plots in Fig.~\ref{fig:s1s2s3MC250SNRn15} for the chosen ASNR and RMSE ranges.   The MVT CRB is shown in Figs.  \ref{fig:s1s2s3MC250SNRn15}(b,e,h).
We have not evaluated the CRB for $\epsilon$-contaminated array data.  In the performance plots for $\epsilon$-contaminated array data,  we show $C_{\mathrm{CR,Gauss}}(\Vec{\theta})$ for $\mathrm{ASNR}=N/\sigma^2$ using \eqref{eq:epscont-variance},  labeled  as  ``Gaussian CRB (shifted)'' in Figs.~\ref{fig:s1s2s3MC250SNRn15}(c,f,i). We expect that the true CRB for $\epsilon$-contaminated array data is higher than this approximation.

\begin{table}
\begin{center}
\normalsize
\begin{tabular}{|l|l|l|} \hline
scenario & DOAs & source variance \\ \hline
single source & $-10^{\circ}$ & $\gamma_{8001}=1$ \\
two sources & $-10^{\circ}$, $10^{\circ}$ & $\gamma_{8001}=\gamma_{10001}=\frac12$\\
three sources & $-3^{\circ}$,  $2^{\circ}$, $75^{\circ}$ & $\gamma_{8701}=\gamma_{9201}=\gamma_{16501}=\frac13$\\ \hline
\end{tabular}
\\*[0.5ex]
\caption{Source scenarios,  source variances normalized to $\mathop{\mathrm{tr}}(\Mat{\Gamma})=1$ } 
\label{tab:source-scenarios}
\end{center}
\end{table}

\begin{table}
\begin{center}
\normalsize
\begin{tabular}{|l|l|l|l|} \hline
array data model & Eq.  & parameters & ASNR \\ \hline
Gaussian & \eqref{eq:linear-model3} & $\Mat{\Sigma}$ & $N/\sigma^2$ \\
MVT &  \eqref{eq:scale-mixture} & $\Mat{\Sigma}$,  $\nu_{\mathrm{data}}=2.1$ & $N/\sigma^2$ \\
$\epsilon$-contaminated & \eqref{eq:linear-model3} & $\Mat{\Sigma}$,  $\epsilon=0.05$, $\lambda=10$ & $N / \sigma^2$ \\ \hline
\end{tabular}
\\*[0.5ex]
\caption{Array data models}
\label{tab:data-models}
\end{center}
\end{table}

\subsection{Single source scenario}
\label{eq:1src}

A single plane wave ($K=1$) with complex circularly symmetric zero-mean Gaussian amplitude is arriving from DOA $\theta_{8001} = -10^{\circ}$.
Here, $\mathrm{ASNR} = N / \sigma^2$, cf. \cite[Eq.~(8.112)]{VanTreesBook}.
Figure \ref{fig:s1s2s3MC250SNRn15} shows results for RMSE of DOA estimates in scenarios with $L=25$ snapshots and $N=20$ sensors.
RMSE is averaged over $250$ iid realizations of DOA estimates from array data $\Mat{Y}$.
There are more snapshots $L$ than sensors $N$, ensuring full rank 
$\Mat{R}_{\Vec{Y}}$ almost surely.
\par
Simulations for Gaussian noise are shown in Fig.~\ref{fig:s1s2s3MC250SNRn15}(a). For this scenario,  the conventional beamformer (not shown) is the ML DOA estimator and approaches the CRB for ASNR greater than $3\,$dB. 
All shown M-estimators for DOA perform almost identically  in terms of RMSE  and just slightly worse than the CBF.
\par 
Figure \ref{fig:s1s2s3MC250SNRn15}(b) shows simulations
for heavy-tailed MVT distributed array data with degrees of freedom parameter $\nu_{\mathrm{data}}=2.1$. 
We observe that the M-estimator for MVT-loss $\rho_{\mathrm{MVT}}$ for loss parameter $\nu_{\mathrm{loss}}=2.1$ 
performs the best, closely followed by the M-estimators for Tyler loss and Huber-loss $\rho_{\mathrm{Huber}}$ 
 for loss parameter $q=0.9$.
Here, the loss parameter $\nu_{\mathrm{loss}}$ used by M-estimator 
is identical to the MVT array data parameter $\nu_{\mathrm{data}}$ and, thus,  is expected to work well. 
In Fig.~\ref{fig:s1s2s3MC250SNRn15}(b), the M-estimator for MVT-loss $\rho_{\mathrm{MVT}}$ closely follows the corresponding MVT CRB \eqref{eq:CES-CRB} for $\mathrm{ASNR}>3\,\mathrm{dB}$, although a small gap at high ASNR remains. 
The assumption that $\nu_{\mathrm{data}}$ is known \emph{a priori} is somewhat unrealistic.
However, methods to estimate $\nu_{\mathrm{data}}$ from data are available,  e.g.,  
\cite{pascal2021improved}.
Gauss loss exhibits largest RMSE at high ASNR in this case.
\par
Results for $\epsilon$-contaminated noise are shown in Fig.~\ref{fig:s1s2s3MC250SNRn15}(c) for outlier probability $\epsilon=0.05$ and outlier strength $\lambda=10$. 
The resulting noise variance \eqref{eq:epscont-variance} for this heavy-tailed distribution is $\sigma^2=5.95\sigma^2_1$. 
The M-estimators for Tyler loss and MVT-loss $\rho_{\mathrm{MVT}}$ perform best followed by
the M-estimator for Huber loss $\rho_{\mathrm{Huber}}$ with an ASNR penalty of about 2\,dB.  
Worst RMSE exhibits the (non-robust) DOA estimator for Gauss loss $\rho_{\mathrm{Gauss}}$ indicating strong impact of outliers on the
original (non-robust) SBL algorithm for DOA estimation \cite{gerstoft2016mmv}. 

\subsection{Two source scenario} \label{section:2src}

Next,  we consider the two-source scenario ($K=2$) in Table \ref{tab:source-scenarios}. 
Here, $\mathrm{ASNR} = N / \sigma^2$, cf. \cite[Eq.~(8.112)]{VanTreesBook}.
Figure \ref{fig:s1s2s3MC250SNRn15} shows results for RMSE of DOA estimates in scenarios with $L=25$ snapshots and $N=20$ sensors. 
\par
Simulations for Gaussian noise are shown in Fig.~\ref{fig:s1s2s3MC250SNRn15}(d). Here, the M-estimate for Huber loss $\rho_{\mathrm{Huber}}$  for loss parameter $q=0.9$ performs equally well as the M-estimate for Gauss loss $\rho_{\mathrm{Gauss}}$ which is equivalent to the original (non-robust) SBL algorithm for DOA estimation \cite{gerstoft2016mmv}. 
They approach the CRB for ASNR greater $9\,$dB.  The DOA estimator for Tyler loss 
performs slightly worse than the previous two. 
Here,   MVT-loss $\rho_{\mathrm{MVT}}$ for loss parameter $\nu_{\mathrm{loss}}=2.1$ has highest RMSE in DOA M-estimates.
\par 
Figure \ref{fig:s1s2s3MC250SNRn15}(e) shows simulations
for heavy-tailed MVT array data with  parameter $\nu_{\mathrm{data}}=2.1$ being small. 
We observe that M-estimation with Tyler loss and MVT-loss $\rho_{\mathrm{MVT}}$ for loss parameter $\nu_{\mathrm{loss}}=2.1$ perform best,  closely followed by M-estimation with Huber loss $\rho_{\mathrm{Huber}}$ with loss parameter $q=0.9$.  The non-robust DOA estimator for $\rho_{\mathrm{Gauss}}$ performs much worse than the other two showing an ASNR penalty of  about 6\,dB at high ASNR.

Here, the loss parameter $\nu_{\mathrm{loss}}$ used by the M-estimator for MVT-loss $\rho_{\mathrm{MVT}}$ is identical to the MVT array data parameter $\nu_{\mathrm{data}}$ used in generating the MVT-distributed array data.
In Fig.~\ref{fig:s1s2s3MC250SNRn15}(e),  the RMSE for  $\rho_{\mathrm{MVT}}$ closely follows the MVT CRB \eqref{eq:CES-CRB} for $\mathrm{ASNR}>9\,\mathrm{dB}$, although a small gap at high ASNR remains.  
\par
Results for $\epsilon$-contaminated noise are shown in Fig.~\ref{fig:s1s2s3MC250SNRn15}(f) for outlier probability $\epsilon=0.05$ and outlier strength $\lambda=10$. 
M-estimation with Tyler loss and MVT-loss $\rho_{\mathrm{MVT}}$ for loss parameter $\nu_{\mathrm{loss}}=2.1$ show lowest RMSE followed by  Huber loss with slight ASNR penalty.  
Gauss loss has again the worst performance.

\subsection{Three source scenario}
\label{eq:3src}

Array data $\Mat{Y}$ are generated for the three source scenario ($K=3$),  so that $\trace{\Mat{\Gamma}}=1$. 
The RMSE performance shown in Figs.  \ref{fig:s1s2s3MC250SNRn15}(g,h,i)  are very similar to the corresponding results for the two-source scenario shown in Figs. \ref{fig:s1s2s3MC250SNRn15}(d,e,f).

\subsection{Effect of Loss Function}

\begin{figure}[t]
\center{\includegraphics[width=.8\columnwidth]{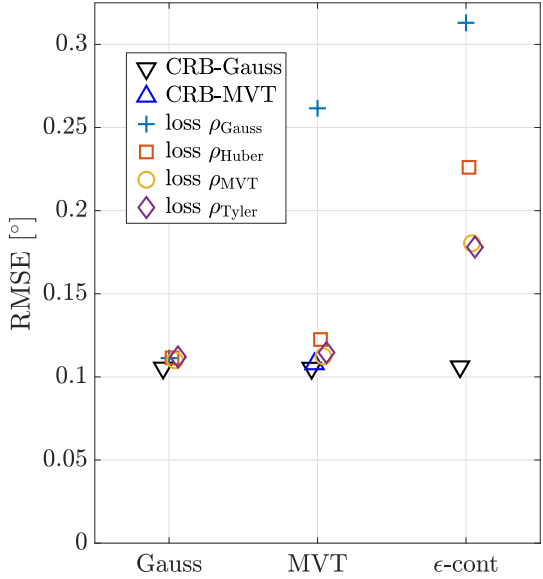}}
       \caption{RMSE for each DOA M-estimator at high $\mathrm{ASNR}=30\,$dB. for each of three array data models (Gaussian,  MVT and $\epsilon$-contaminated).}
\label{fig:summary}
\end{figure}

The effect of the loss function on RMSE performance at high $\mathrm{ASNR}=30\,$dB is illustrated in Fig. ~\ref{fig:summary}. 
This shows that for Gaussian array data all choices of loss functions perform equally well at high ASNR.   For MVT data in Fig.~\ref{fig:summary} (middle), we see that the robust loss functions (MVT,  Huber,  Tyler) work well,  and approximately equally,  whereas RMSE for Gauss loss is factor 2 worse. 
For $\epsilon$-contaminated array data in Fig.~\ref{fig:summary} (right) the Gauss loss performs a factor worse than the robust loss functions.  Huber loss has slightly higher RMSE than MVT and Tyler loss.

\subsection{Effect of Outlier Strength on RMSE}
\label{sec:outlier-strength}

\begin{figure}[t]
{\footnotesize \textcolor{white}{\tt RMSEvsLambda.eps}}\\
\includegraphics[width=\columnwidth]{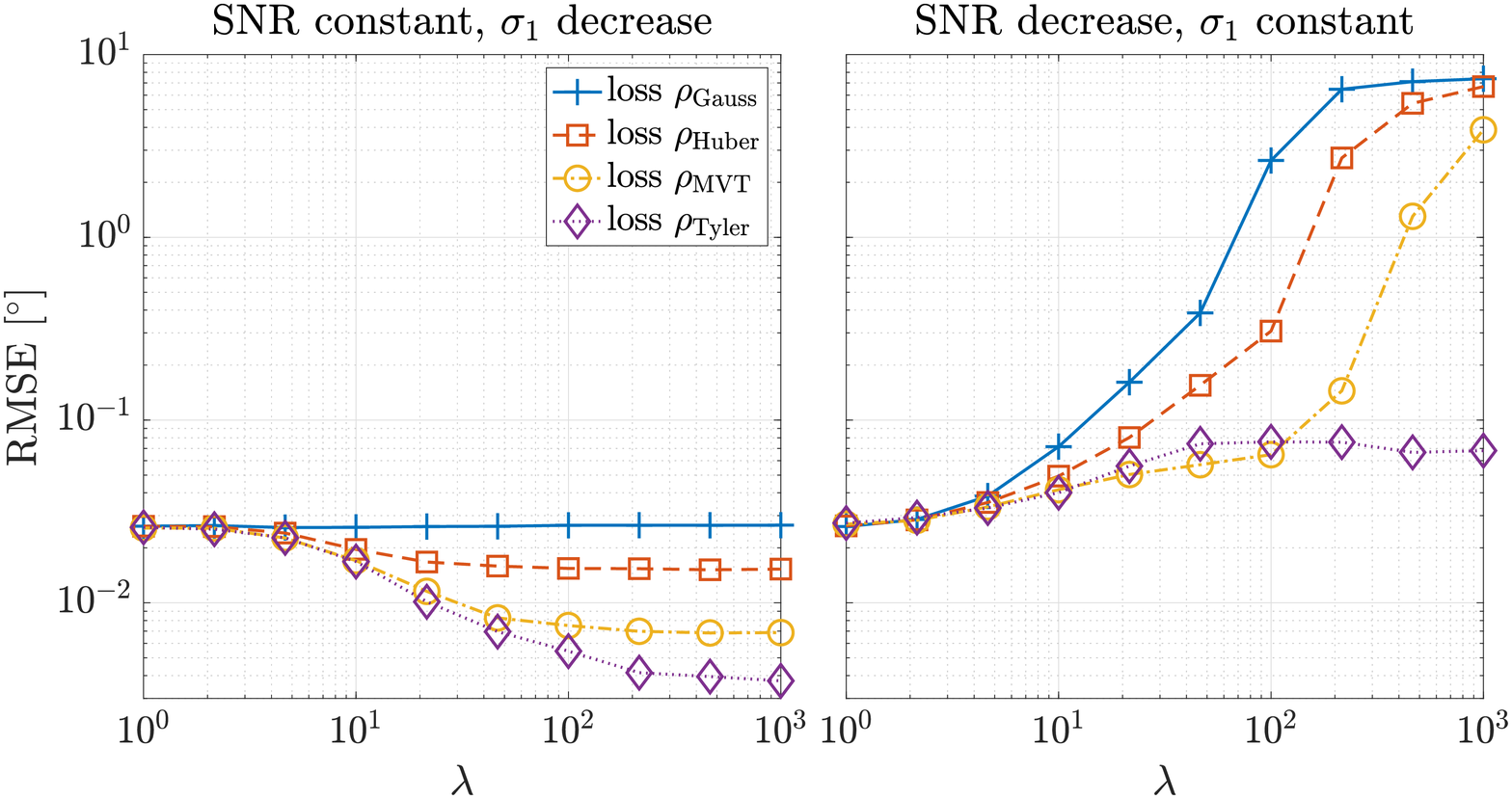}
\caption{\label{fig:RMSEvsLambda}  For $\epsilon$-contaminated array  data $\epsilon=0.05$ with one source, RMSE versus outlier strength $\lambda$ for each loss function for (left): ASNR$= 25\,$dB and background noise $\sigma_1$ is decreasing, 
and (right): Background noise $\sigma_1$ is fixed and outlier noise  $\lambda\sigma_1$ is increasing (at $\lambda=1$: ASNR$=25\,$dB;  and at $\lambda=10^3$: ASNR $=25-10 \log(1-\epsilon+\epsilon\lambda^2)=-22$ 
dB ).  RMSE evaluation based on $N_{\mathrm{run}}=250$ simulation runs.}
\end{figure}

For small outlier strength $\lambda$ for $\epsilon$-contaminated data, the Gauss loss performs satisfactorily,  but as the outlier noise increases the robust loss functions  outperform, see Fig.~\ref{fig:RMSEvsLambda}.
As $\lambda$  increases,  the total noise changes,  see \eqref{eq:epscont-variance}. 
In Fig.~\ref{fig:RMSEvsLambda}(left) the total noise is kept constant by decreasing the background noise with increasing   $\lambda$.   In Fig.~\ref{fig:RMSEvsLambda}(right) the background noise level is constant  whereby the total noise increases.
For large noise outlier,  Tyler loss clearly has best performance in Fig.~\ref{fig:RMSEvsLambda}(left) and does not breakdown in Fig.~\ref{fig:RMSEvsLambda}(right).

\begin{figure}[t]
{\footnotesize \textcolor{white}{\tt OtherAlg.eps}}\\
\includegraphics[width=\columnwidth]{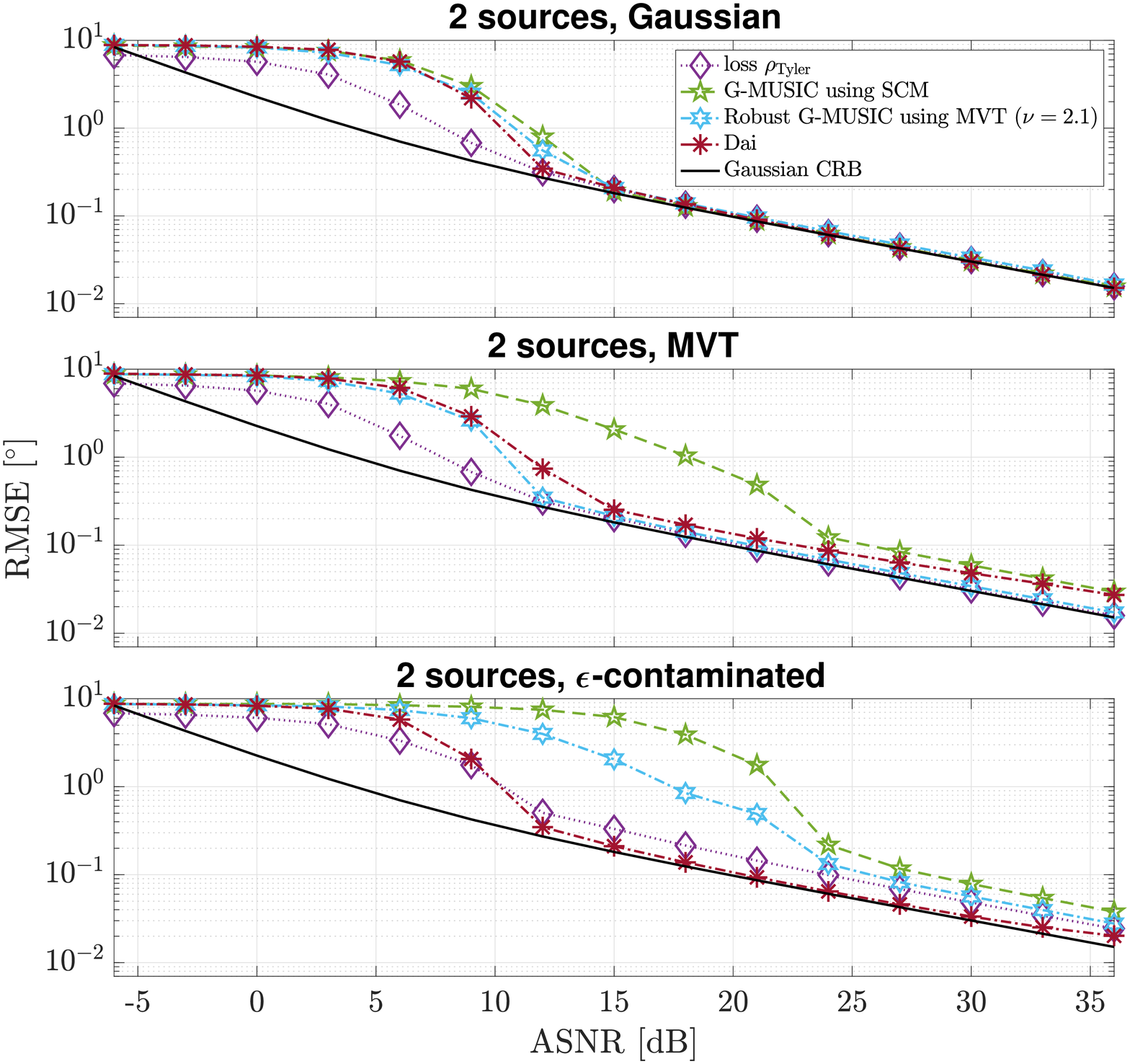}
\caption{\label{fig:otheralg} RMSE vs ASNR for Gaussian, MVT, and $\epsilon$-contaminated array data ($\epsilon=0.05, \lambda=10$).  Similar setup as Fig.~\ref{fig:s1s2s3MC250SNRn15}  with 2 sources at $\pm10^\circ$. }
\end{figure}

\subsection{Comparison with other algorithms}
\begin{figure}[t]
{\footnotesize \textcolor{white}{\tt CORRsrc.eps}}\\
\includegraphics[width=\columnwidth]{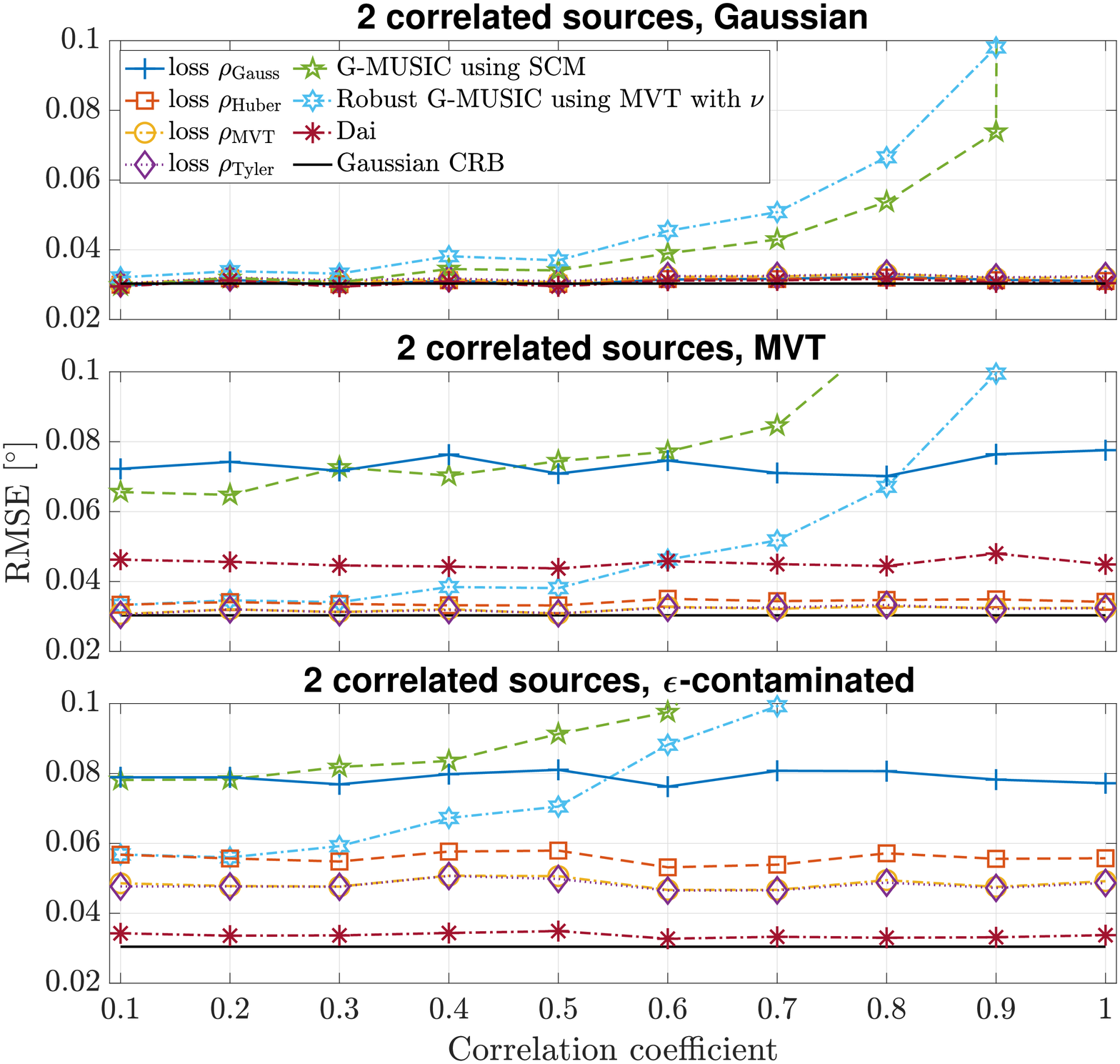}
\caption{\label{fig:coherent} RMSE vs source correlation for Gaussian, MVT, and $\epsilon$-contaminated array data ($\epsilon=0.05, \lambda=10$) at ASNR$=30\,$dB.  Similar setup as Fig.\ \ref{fig:s1s2s3MC250SNRn15}  with 2 sources at $\pm10^\circ$. }
\end{figure}

The RMSE of DOA estimates from Table \ref{table:algorithm} is compared with G-MUSIC with SCM  \eqref{eq:SY} \cite{Mestre2008}, G-MUSIC \cite{Couillet2014}, as well as with Dai and So \cite{DaiSo2018} in Figs.
\ref{fig:otheralg} and \ref{fig:coherent}.
The Dai and So algorithm \cite{DaiSo2018} performs well for $\epsilon$-contaminated noise,  as the outliers are assumed Gaussian.

\subsection{Coherent sources}
Coherent sources are often encountered in applications,  e.g.,  
ocean acoustics and radar array processing.
In contrast to subspace methods 
SBL performs well for coherent sources \cite{das2017, mecklenbrauker2018,pote2020robustness}.
To demonstrate this,  we report results for a two-source scenario with DOAs at $\pm 10^\circ$ with source correlation coefficient 
\begin{equation}
\rho = \frac{ \mathsf{E}[x_{8001,l}x_{10001,l}^*]}{\sqrt{\mathsf{E}[ |x_{8001,l}|^2] \, \mathsf{E}[ |x_{10001,l}|^2]}}.
\end{equation}
The data model for coherent sources has a full source correlation matrix $\Mat{\Gamma}$, but it has been found that estimating a full matrix $\Mat{\Gamma}$ is numerically costly
 \cite{mecklenbrauker2018,pote2020robustness} and not actually needed for DOA estimation.
The Robust SBL algorithm in Table  \ref{table:algorithm} processes coherent sources by ignoring any off-diagonal elements in $\Mat{\Gamma}$, i.e.,  by using the diagonal $\Mat{\Gamma}$ model \eqref{eq:Gamma-model_new}. 

The numerical results in Fig.~\ref{fig:coherent} demonstrate that the coherency of the sources has little effect on the SBL results from Table \ref{table:algorithm} and this is also true for Dai and So \cite{DaiSo2018},  in contrast to G-MUSIC with either SCM  \eqref{eq:SY}  \cite{Mestre2008} or \cite{Couillet2014} when $\rho>0.5$.

\subsection{Effect of Dictionary Size on RMSE}
\label{sec:hires-grid}
Due to the algorithmic speedup associated with the DOA grid pruning described in 
Sec.  \ref{sec:algorithm},  it is feasible and useful to run the algorithm in Table \ref{table:algorithm} with large dictionary size $M$
which translates to the dictionary's angular grid resolution of $\delta=180^\circ/(M-1)$.

\begin{figure}[t]
{\footnotesize \textcolor{white}{\tt gridresolution.eps}}\\
\includegraphics[width=\columnwidth]{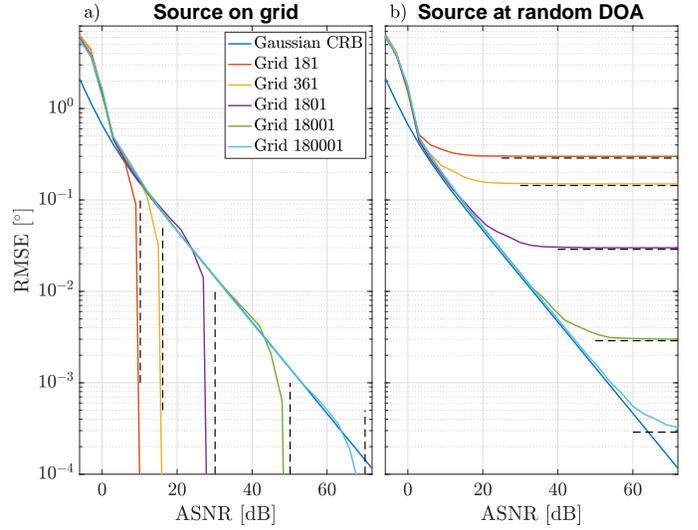}
\caption{\label{fig:gridresolution} Effect of dictionary size $M\in\{181, 361, 1801, 18001, 180001 \}$ on RMSE vs.~SNR for a single source a) fixed DOA $-10^\circ$ on the grid and b) random uniformly distributed DOA $\sim -10^\circ + U(-\delta/2,\delta/2)$.
The vertical and horizontal dashed lines indicate asymptotic values due to dictionary size.}
\end{figure}

The effect of grid resolution is illustrated in Fig.~\ref{fig:gridresolution} for a single source impinging on a $N=20$ element $\lambda/2$-spaced ULA.
The Gaussian array data model is used. 
Fig.~\ref{fig:gridresolution} shows RMSE vs.\  ASNR for a dictionary size  of $M\in\{181, 361, 1801, 18001, 180001 \}$. 
In Fig.~\ref{fig:gridresolution}(a), the DOA is fixed
at $-10^\circ$, cf.  single source scenario in Table \ref{tab:source-scenarios}, the DOA is on the angular grid which defines the dictionary matrix $\Mat{A}$.
In Fig.~\ref{fig:gridresolution}(b) 
the  DOA is random,  the source DOA  is sampled from $-10^\circ + U(-\delta/2,  \delta/2)$
($\delta=180^{\circ}/(M-1)$ is the angular grid resolution).
The source DOA is not on the angular grid which defines the dictionary matrix $\Mat{A}$.

For source DOA on the dictionary grid,  Fig.~\ref{fig:gridresolution}(a),    the RMSE performance curve resembles the behavior of an ML-estimator at low ASNR up to a certain threshold ASNR (dashed vertical lines) where the RMSE abruptly crosses the CRB and becomes zero. 
The threshold ASNR is deduced from the following argument:
Let $\Vec{a}_m$ be the true DOA dictionary vector and $\Vec{a}_{m+1}$ be the dictionary vector for adjacent DOA on the angular grid.  Comparing the corresponding Bartlett powers,  we see that DOA errors become likely if the noise variance exceeds $2(|\Vec{a}_m^{\sf H}\Vec{a}_m| -|\Vec{a}_m^{\sf H} \Vec{a}_{m +1}|)/N = 2 - 2|\Vec{a}_m^{\sf H}\Vec{a}_{m +1}|/N$.

For source DOA off the dictionary grid,  Fig.~\ref{fig:gridresolution}(b),   the RMSE performance curve resembles the behavior of an ML-estimator at low ASNR up to a threshold ASNR.  In the random DOA scenario,  however,  the RMSE flattens at increasing ASNR. 
Since the variance of the uniformly distributed source DOA is $\delta^2/12$, 
the limiting $\mathrm{RMSE}=\delta/\sqrt{12}$ for $\mathrm{ASNR}\to\infty$.
The limiting RMSE (dashed horizontal lines) depends on the dictionary size $M$ through the angular grid resolution $\delta$.   
The asymptotic RMSE limits are shown as dashed horizontal lines in Fig.~\ref{fig:gridresolution}(b).

\section{Convergence Behavior and Run Time}
\label{sec:convergence}

The DOA M-estimation algorithm in Table \ref{table:algorithm} uses an iteration to estimate the active set $\mathcal{M}$ whose elements represent the estimated source DOAs.  The required number of iterations for convergence of $\mathcal{M}$
 depends on the source scenario,  array data model,  ASNR,   and stepsize $\mu$.  We select $\mu=1$.

\begin{figure}[t]
{\footnotesize \textcolor{white}{\tt noItr.eps}}\\
\includegraphics[width=\columnwidth]{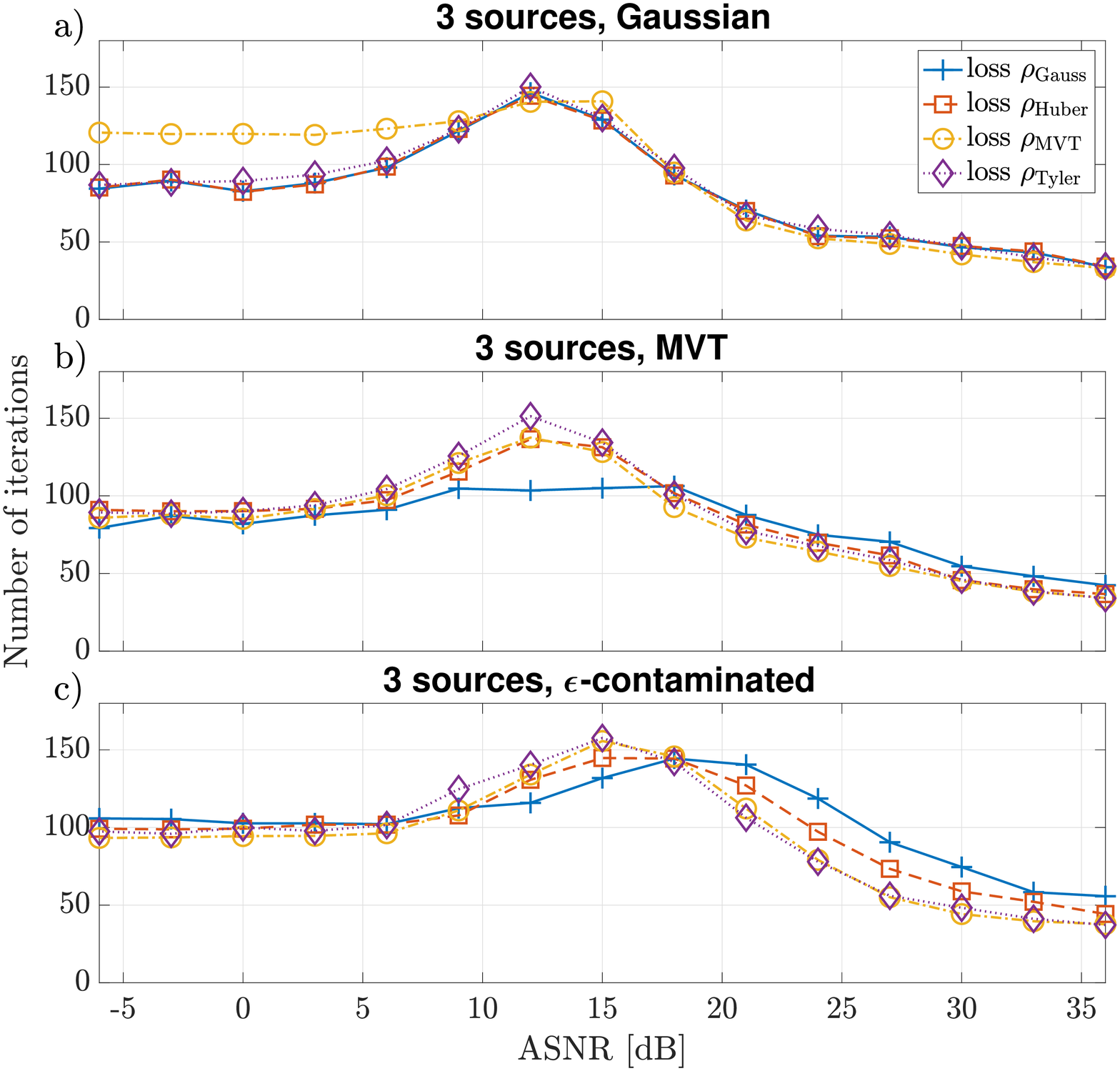}
\caption{\label{fig:s3MC250SNRn15cpu} Iteration count of DOA estimators vs.~ASNR for three sources,  $N=20$ sensors, $L=25$ array snapshots, and dictionary size $M=18001$ corresponding to DOA resolution $0.01^{\circ}$. 
Noise: (a) Gaussian, (b) MVT ($\nu_{\mathrm{data}}=2.1$), (c) $\epsilon$-contaminated   ($\epsilon=0.05,\lambda=10$).}
\end{figure}

Figure \ref{fig:s3MC250SNRn15cpu} shows the required number of iterations for the 
three-source scenario versus ASNR and all three array data models.  
Figure \ref{fig:s3MC250SNRn15cpu}  shows fast convergence for high ASNR, and the number of iterations decreases with increasing ASNR.  
At $\mathrm{ASNR}< 5\,$dB,  where the noise dominates,  the number of iterations is around 100 and approximately indepent of the ASNR.
Figure \ref{fig:s3MC250SNRn15cpu}(a) shows that  the number of iterations for MVT-loss at low ASNR and Gaussian array data is about 25\% larger than for the other loss functions.
In the intermediate ASNR range 5--20\,dB,  the largest number of iterations are required as the algorithm searches the dictionary to find the best matching DOAs.  
Peak number of iterations is near 160 at ASNR levels between 12 and 15 dB. 
Figure \ref{fig:s3MC250SNRn15cpu}(c) shows that  the number of iterations for Tyler loss and MVT-loss for $\epsilon$-contaminated array data
at high ASNR are lowest, followed by Huber loss and Gauss loss.

\begin{figure}[t]
{\footnotesize \textcolor{white}{\tt CPUtime.eps}}\\
\includegraphics[width=\columnwidth]{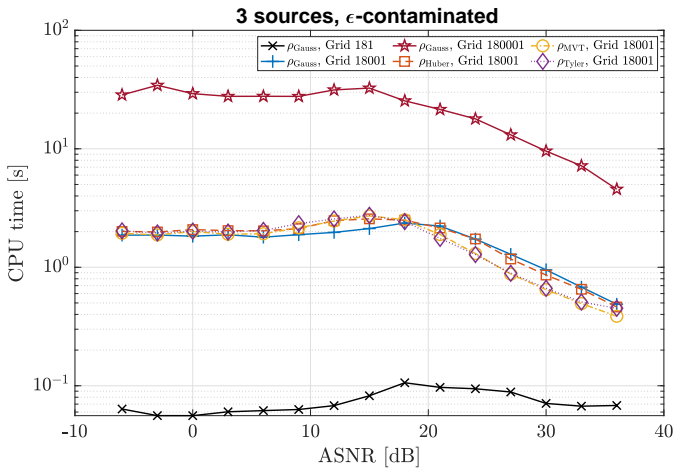}
\caption{\label{fig:CPUtimes} CPU times for three source scenario vs.  ASNR for $\epsilon$-contaminated array data processed with Gauss loss,  dictionary size $M\in\{181, 18001, 180001 \}$ and for  $M=18001$ for Huber, MVT and Tyler loss.}
\end{figure}

The CPU times on an M1 MacBook Pro are shown in Fig.~\ref{fig:CPUtimes} for $\epsilon$-contaminated array data and various choices of dictionary size and loss function. 
For a fixed dictionary ($M=18001$), the choice of loss function does not have much influence on the CPU time.  
At $\mathrm{ASNR}> 18\,$dB,  MVT and Tyler consume just slightly less CPU time than Huber and Gauss loss.
For low ASNR,  CPU time increases by a ratio that is approximately proportional to dictionary size 180001/181=1000, but at high SNR this ratio reduces to 50, due to the efficiency of the DOA grid pruning,  cf.  Sec.  \ref{sec:algorithm}.
The $M=180001$ dictionary is quite large in this scenario,  but in other scenarios for localization in 3 dimensions,  this is an expected dictionary size.

\section{Conclusion}
\label{sec:conclusion}

Robust and sparse DOA M-estimation is derived based on array data following a zero-mean complex elliptically symmetric  distribution with finite second-order moments.   Our derivations allowed using different loss functions. 
The DOA M-estimator is numerically evaluated by iterations 
and available as a Matlab function on GitHub \cite{RobustSBL-github}.
A specific choice of loss function determines the RMSE performance and robustness of the DOA estimate. 
Four choices for loss function are discussed and investigated in numerical simulations with synthetic array data: the ML-loss function for the circular complex multivariate $t$-distribution with $\nu$ degrees of freedom,  the loss functions for Huber and Tyler M-estimators.  For Gauss loss, the method reduces to Sparse Bayesian Learning.
We discuss the robustness of these  DOA M-estimators by evaluating the root mean square error for Gaussian,  MVT,  and $\epsilon$-contaminated array data.  The robust and sparse M-estimators for DOA perform well in simulations for MVT and $\epsilon$-contaminated noise and nearly identical with classical SBL for Gaussian noise. 

\section*{Acknowledgement}
The authors thank Markus Rupp for sharing valuable insights on convergence and the anonymous reviewers for helping to improve the paper.

\appendix
\subsection{MVT array data model}
\label{sec:MVT-array-data-model}

Here we show how the MVT array data model relates to the 
scale mixture model for the array data \eqref{eq:scale-mixture}.
By assuming that the distribution of $\tau_{\ell}$  in \eqref{eq:scale-mixture}
is an inverse gamma distribution with shape parameter $\alpha$ and scale parameter $\beta$,
\begin{align}
   \tau_{\ell} &\sim p_{\tau}(\tau) = \frac{\beta^\alpha}{\Gamma(\alpha)} \frac{1}{\tau^{\alpha+1}} \Exp{-\beta/\tau},
               \text{ for } \tau>0,
\end{align}
the corresponding density generator evaluates to
 \begin{align} \label{eq:t_g2}
 g(t) &= \pi^{-N} \int_0^\infty \tau^{-N} e^{-t/\tau} p_{\tau}(\tau)  \,\mathrm{d}\tau \\
        &=  \frac{\beta^\alpha}{\pi^N\Gamma(\alpha)}   
               \int_0^\infty  \frac{e^{-(t+\beta)/\tau}}{\tau^{N+\alpha+1}} \,\mathrm{d}\tau.
\end{align}
We substitute $w = 1/\tau$ and $\mathrm{d}\tau = - \frac{1}{w^2}\mathrm{d}w$, giving
 \begin{align} \label{eq:t_g3}
 g(t)  &=  \frac{\beta^\alpha}{\pi^N\Gamma(\alpha)}   
               \int_0^\infty  w^{N+\alpha-1} \Exp{-(t+\beta)w} \,\mathrm{d}w \\
 &=  \frac{\beta^\alpha}{\pi^N\Gamma(\alpha)}   
             \frac{\Gamma(N+\alpha)}{(t+\beta)^{N+\alpha}}.
\end{align}
Finally,  we specialize the shape and scale parameters to $\alpha=\beta=\nu_{\mathrm{data}}/2$, giving
\begin{align} \label{eq:t_g4}
 g(t) &=  \frac{(\nu_{\mathrm{data}}/2)^{\nu_{\mathrm{data}}/2}}{\pi^N\Gamma(\nu_{\mathrm{data}}/2)}   
             \frac{\Gamma(N+\nu_{\mathrm{data}}/2)}{(\nu_{\mathrm{data}}/2+t)^{N+\nu_{\mathrm{data}}/2}} \\
      &=  \left( \frac{2}{\pi\nu_{\mathrm{data}}} \right)^N
            \frac{\Gamma(N+\nu_{\mathrm{data}}/2)}{\Gamma(\nu_{\mathrm{data}}/2)}   
             \left(1+\frac{2t}{\nu_{\mathrm{data}}}\right)^{-N-\nu_{\mathrm{data}}/2} 
\end{align}
which agrees with the density generator of the multivariate $t$-distribution  $t_{\nu_{\mathrm{data}}}(\Vec{0},\Vec{\Sigma})$ \cite[Sec. 4.2.2,  Ex.  11 on p.  107]{Zoubir2018}.  Thus, placing an inverse gamma prior over the random scaling $\tau_{\ell}$ results in the MVT array data model.

\subsection{Consistency Factor}
\label{sec:consistency-factor}

Here the consistency factor $b$ is evaluated for  Huber loss  \eqref{eq:loss-H},  MVT loss, \eqref{eq:loss-T},  and Tyler loss \eqref{eq:loss-Tyler}.  For elliptical distributions an M-estimator is a consistent estimator of $\alpha \Mat{\Sigma}$, 
where the constant $\alpha $ is a solution to \cite[eq. (49)]{ollila2012complex}:
\begin{equation} \label{eq:b_estim_equation}
1 = \mathsf{E}  \big[ \psi(  \Vec{y}^{\sf H} \Mat{\Sigma}^{-1} \Vec{y} / \alpha) \big]/N, 
\end{equation} 
where $\psi(t)= t \,u(t) = t \,\mathrm{d}\rho(t)/\mathrm{d}t$ as defined below \eqref{eq:consitency factor2} .  

Assuming $\Vec{y} \sim \mathbb{C} \mathcal N_N(\Vec{0},\Mat{\Sigma})$,  then we scale the chosen loss function $\rho(t)$ such that 
\eqref{eq:b_estim_equation} holds for $\alpha=1$. Namely,  for 
\begin{equation}
 \rho_b(t) = \rho(t)/b  \mbox{ and } u_b(t) = u(t)/b, 
\end{equation}
where $b$ is a scaling constant defined in \eqref{eq:consitency factor}, it clearly holds that 
$1 = \mathsf{E}  \big[ \psi_b(  \Vec{y}^{\sf H} \Mat{\Sigma}^{-1} \Vec{y}) \big]/N$ for  
$\psi_b(t)= t u_b(t) $.  This implies that $\alpha=1$ and that the M-estimator  with loss $\rho_b(\cdot)$ is consistent to the covariance matrix $\Mat{\Sigma}$ when the array data follows the  $\mathbb{C} \mathcal N_N(\Vec{0},\Mat{\Sigma})$ distribution. 
 
For Huber loss function \eqref{eq:loss-H} the $b$  in \eqref{eq:consitency factor2}  can be solved in closed form as \cite[Sec. 4.4.2]{Zoubir2018}
\begin{align} 
b &=  \frac 1 N  \int_0^\infty  (t/2) u_{\mathrm{Huber}}( t/2)  f_{\chi^2_{2N}}(t) \mathrm{d} t   \\ 
&=  \frac 1 {2 N}   \int_0^{2 c^2}  t f_{\chi^2_{2N}}(t) \mathrm{d} t   +   \frac 1 { N}   \int_{2 c^2}^{\infty} c^2  f_{\chi^2_{2N}}(t) \mathrm{d} t   \\ 
&=  F_{\chi^2_{2(N+1)}}(2c^2) + c^2(1- F_{\chi^2_{2 N}}(2 c^2))/N,
\end{align} 
where $F_{\chi^2_{n}}(x)=\mathrm{P}\{ X \le x \}$ is the cumulative distribution of a central $\chi^2$ distributed random variable $X$ with $n$ degrees of freedom.

For the MVT loss \eqref{eq:loss-T} we evaluate \eqref{eq:consitency factor2}  by numerical integration.  

For Tyler loss $\psi(t)\equiv N$ $\forall t$,  indicating that the consistency factor for Tyler loss can not be found based on \eqref{eq:b_estim_equation}. 
However it is possible to construct an affine equivariant version of Tyler’s M-estimator as explained in \cite{ollila2023affine} that estimates the true scatter matrix $\Mat{\Sigma}$.
This is accomplished by defining $b$ as the mean of Tyler’s weights
\begin{align}
  b_{\mathrm{Tyler}}  &= \frac1L \sum\limits_{\ell=1}^L u_{\mathrm{Tyler}} (\Vec{y}_\ell \Mat{\Sigma}^{-1} \Vec{y}_\ell).
\label{eq:consistency-factor-Tyler}
\end{align}

\subsection{Convergence of the iteration for gamma} 
\label{sec:stability} 
Here, we discuss the convergence of iterations based on \eqref{eq:gamma-update-rule-with-stepsize}.
We proceed in two steps:

{\it 1. True solution $\Vec{\gamma}^{\mathrm{true}}$ is a fixed point of the iteration \eqref{eq:gamma-update-rule-with-stepsize}:}
It follows from \eqref{eq:gamma-update-rule-with-stepsize} that
$\Vec{\gamma}^{\mathrm{new}} =  \Vec{\gamma}^{\mathrm{old}}$
when $\Vec{\gamma}^{\mathrm{old}}=\Vec{\gamma}^{\mathrm{true}}$ equals the true vector of source powers $\Vec{\gamma}^{\mathrm{true}}$
 if $\Mat{R}_{\Mat{Y}}=\Mat{\Sigma}$ for any choice of $\mu$. 
 The estimate $\Mat{R}_{\Mat{Y}}$ is unbiased for $\Mat{\Sigma}$ 
 when $\Vec{y}_{\ell}$ are Gaussian, or follow any of distributions in Table \ref{tab:loss-and-weight-functions} thanks to the consistency factor $b$.  
 Asuming $\Mat{R}_{\Mat{Y}}$ is unbiased for $\Mat{\Sigma}$,
\begin{align}
      \mathsf{E} [ G_m(\Vec{\gamma}^{\mathrm{true}}) ] = 1,
\quad\forall m.
\end{align}
Next, we abbreviate $\Vec{b}_m = \Mat{\Sigma}^{-1}\Vec{a}_m$ and rewrite \eqref{eq:iteration-gain-Gm},
\begin{align}
      G_m(\Vec{\gamma}) &=    \frac
       {\Vec{b}_m^{\sf H} \Mat{R}_{\Mat{Y}} \Vec{b}_m}
       {\Vec{b}_m^{\sf H} \, \Mat{\Sigma} \, \Vec{b}_m} \label{eq:bRb/bSb}\\
      &= G_m(\Vec{\gamma}) - G_m(\Vec{\gamma}^{\mathrm{true}}) + 1, \label{eq:Gm-Gm+1}
\end{align}

{\it 2. The iteration update  \eqref{eq:gamma-update-rule-with-stepsize} decreases the error:}
Using \eqref{eq:Gm-Gm+1} 
and \eqref{eq:gamma-update-rule-with-stepsize},  we obtain
\begin{align}
       \gamma_m^{\mathrm{true}} - &\gamma_m^{\mathrm{new}} =     \gamma_m^{\mathrm{true}} - 
         (1-\mu)   \gamma_m^{\mathrm{old}} - \mu  \gamma_m^{\mathrm{old}} G_m(\Vec{\gamma}^{\mathrm{old}}),   \nonumber \\
    &= \! ( \!\gamma_m^{\mathrm{true}} \!  -  \! \gamma_m^{\mathrm{old}} \! ) \! \!
\left[  \! 1 \! + \! \mu
          \frac{G_m \! (\Vec{\gamma}^{\mathrm{true}} \! ) \!  -  \! G_m \! (\Vec{\gamma}^{\mathrm{old}} \! )}{\gamma_m^{\mathrm{true}} - \gamma_m^{\mathrm{old}} } \gamma_m^{\mathrm{old}}  \! \right]  \! \! .
\end{align}
Next,  the mean value theorem of differential calculus is used to express the difference quotient
\begin{align}
    \frac{G_m(\Vec{\gamma}^{\mathrm{true}}) - G_m(\Vec{\gamma}^{\mathrm{old}})}{\gamma_m^{\mathrm{true}}- \gamma_m^{\mathrm{old}} } &=
 \left(\frac{\partial G_m(\Vec{\gamma})}{\partial \gamma_{m}} \right)_{\hspace*{-0.7ex}\Vec{\gamma}=\tilde{\Vec{\gamma}}}
\end{align}
where $\tilde{\Vec{\gamma}}=\Vec{\gamma}^{\mathrm{true}} + \vartheta \gamma_m^{\mathrm{old}}\Vec{\mathbf{e}}_m$,  $0\le\vartheta\le1$,  and $\Vec{\mathbf{e}}_m$ is the $m$th standard basis vector.
For convergence,   the new error $| \gamma_m^{\mathrm{true}} - \gamma_m^{\mathrm{new}}| $ (after updating) must be less than the previous error $| \gamma_m^{\mathrm{true}} - \gamma_m^{\mathrm{old}}| $.
This is ensured if
\begin{align}
\left| 1 + \mu \gamma_m^{\mathrm{old}}
          \left(\frac{\partial G_m(\Vec{\gamma})}{\partial \gamma_{m}} \right)_{\hspace*{-0.7ex}\Vec{\gamma}=\tilde{\Vec{\gamma}}} \right| & < 1,
\end{align}
giving
\begin{align}
0 & \le -\mu\gamma_m^{\mathrm{old}} 
         \left(\frac{\partial G_m(\Vec{\gamma})}{\partial \gamma_{m}} \right)_{\hspace*{-0.7ex}\Vec{\gamma}=\tilde{\Vec{\gamma}}}  < 2.
\end{align}
Therefore,  in a neighborhood of $\Vec{\gamma}^{\mathrm{true}}$,  we  need
\begin{align}
0 & \le
          - \mu \gamma_m\frac{\partial G_m(\Vec{\gamma})}{\partial \gamma_{m}}   < 2. \label{eq:we-check-the-values}
\end{align}
We distinguish two cases for noise free data,  namely $m\not\in\mathcal{M}$ versus $m\in\mathcal{M}$.
In the first case,  the indices are not active, i.e.,  $\gamma_m^{\mathrm{true}}=0$ due to the sparsity model  \eqref{eq:mathcal-M}.  Consequently \eqref{eq:we-check-the-values} is fulfilled 
 for all $m\not\in\mathcal{M}$ and the associated $\gamma_m$ estimates will correctly converge to zero for any 
chosen stepsize $\mu$ with $0<\mu\le 1$.  
We conclude that the set of active indices $\mathcal{M}$ will be correctly estimated for any $\mu$
The $\mathcal{M}$ can be estimated correctly,  even though the $\gamma_m$ estimates 
might not converge to the true source powers for the active indices, $m\in\mathcal{M}$,  when $\mu$ is chosen large.
Estimation of  $\mathcal{M}$ (rather than $\Vec{\gamma}$) is the primary interest,  because the elements of $\mathcal{M}$ correspond to DOAs.

In the second case,  the indices are active, i.e., $\gamma_m>0$ for $m\in\mathcal{M}$.
Using \eqref{eq:bRb/bSb},  we evaluate
\begin{align}
- \frac{\partial G_m(\Vec{\gamma})}{\partial \gamma_{m}} 
&= \frac{(\Vec{b}_m^{\sf H} \Mat{\Sigma}\Vec{b}_m)'(\Vec{b}_m^{\sf H} \Mat{R}_{\Mat{Y}} \Vec{b}_m)}{(\Vec{b}_m^{\sf H} \Mat{\Sigma}\Vec{b}_m)^2} -\frac{(\Vec{b}_m^{\sf H} \Mat{R}_{\Mat{Y}} \Vec{b}_m)'}{\Vec{b}_m^{\sf H} \Mat{\Sigma}\Vec{b}_m} 
\nonumber \\
&= 
     \Vec{b}_m^{\sf H} \Mat{R}_{\Mat{Y}} \Vec{b}_m
        -\frac{\Vec{b}_m^{\sf H} \Mat{R'_Y}\Vec{b}_m}{\Vec{b}_m^{\sf H} \Mat{\Sigma}\Vec{b}_m} 
\label{eq:Gm-derivative}\\
\text{with } (\cdot)' &:= \partial (\cdot) / \partial \gamma_m \nonumber \\
\Vec{b}'_m& = \partial \Vec{b}_m/\partial\gamma_m = 
- \Vec{b}_m (\Vec{b}_m^{\sf H} \Mat{\Sigma}\Vec{b}_m),\\
\Mat{R}'_{\Mat{Y}} &=\partial \Mat{R}_{\Mat{Y}}/\partial\gamma_m\\
\Mat{\Sigma}' &= \Vec{a}_m \Vec{a}_m^{\sf H} \\
(\Mat{\Sigma}^{-1})' &= -\Mat{\Sigma}^{-1} \Mat{\Sigma}' \Mat{\Sigma}^{-1} = -\Vec{b}_m \Vec{b}_m^{\sf H} \\
(\Vec{b}_m^{\sf H} \Mat{\Sigma}\Vec{b}_m)'&=-(\Vec{b}_m^{\sf H}\Mat{\Sigma}\Vec{b}_m)^2 \\
(\Vec{b}_m^{\sf H} \Mat{R}_{\Mat{Y}} \Vec{b}_m)' 
&=\Vec{b}_m^{\sf H} \Mat{R}_{\Mat{Y}} '\Vec{b}_m    -2\Vec{b}_m^{\sf H} \Mat{R}_{\Mat{Y}} \Vec{b}_m(\Vec{b}_m^{\sf H} \Mat{\Sigma}\Vec{b}_m)
\end{align}
For evaluating the second term  in \eqref{eq:Gm-derivative}, we use \eqref{eq:RY}, giving
\begin{align}
\Mat{R'_Y} &= \frac{1}{L} \Mat{Y}\Mat{D}' \Mat{Y}^{\sf H},  \text{with } \Mat{D}' = \mathop{\mathrm{diag}}(d_1,\ldots,d_L)/b,
\\
\text{with } d_\ell &= u'(\Vec{y}_\ell^{\sf H} \Mat{\Sigma}^{-1} \Vec{y}_\ell) 
(\Vec{y}_\ell^{\sf H} (\Mat{\Sigma}^{-1})'  \Vec{y}_\ell)
\nonumber \\
   &= - u'(\Vec{y}_\ell^{\sf H} \Mat{\Sigma}^{-1} \Vec{y}_\ell) \, (\Vec{y}_\ell^{\sf H} \Vec{b}_m\Vec{b}_m^{\sf H} \Vec{y}_\ell) \nonumber\\
   &= - u'(\Vec{y}_\ell^{\sf H} \Mat{\Sigma}^{-1} \Vec{y}_\ell) \, |\Vec{y}_\ell^{\sf H} \Vec{b}_m|^2 
\end{align}
For Gauss loss,  $\Mat{R}'_{\Mat{Y}} = \Mat{0}$.  
Since all weight functions in Table \ref{tab:loss-and-weight-functions} can be written in the form $u(t)=1/(a+bt)$ with $a,b\ge0$,
it follows that $-u'(t)$ is non-negative and bounded for $t\ge0$.
We conclude that $\Mat{R}'_{\Mat{Y}}$ is positive semi-definite and second term  in \eqref{eq:Gm-derivative} is non-negative,
\begin{align}
\frac{\Vec{b}_m^{\sf H} \Mat{R}'_{\Mat{Y}} \Vec{b}_m}{\Vec{b}_m^{\sf H} \Mat{\Sigma}\Vec{b}_m} &\ge 0.
\end{align}
Continuing from \eqref{eq:Gm-derivative}, we now use the lower bound
\begin{align}
- \frac{\partial G_m(\Vec{\gamma})}{\partial \gamma_{m}}\gamma_m   &\ge \max\limits_{m\in\mathcal{M}} \left(  {\Vec{b}_m^{\sf H}  \Mat{T}_{\Mat{Y}} \Vec{b}_m}
                        \gamma_m \right) = \beta, \\
\text{with } \Mat{T}_{\Mat{Y}} &= \Mat{R}_{\Mat{Y}} - \frac{\Mat{R}'_{\Mat{Y}}}{\Vec{b}_m^{\sf H} \Mat{\Sigma}\Vec{b}_m},
\end{align}
and $\beta >0$ is needed.  This is ensured when $\Mat{T}_{\Mat{Y}}$ is positive semi-definite. 
Returning to \eqref{eq:we-check-the-values}, we see that any chosen stepsize
\begin{align}
      0 &< \mu <
\frac{2}{\beta} 
\label{eq:convergence-condition} 
\end{align}
guarantees convergence for $m\in\mathcal{M}$.

Finally,  we evaluate $\mathsf{E}[\beta]$ for a source scenario with $K$ incoherent sources, 
$\mathcal{M}=\{1,\ldots,K\}$,
$\Mat{\Gamma}_\mathcal{M} =\mathop{\mathrm{diag}}[\gamma_1, \ldots, \gamma_K]$,
assume $\mathsf{E}[\Mat{R}_{\Mat{Y}}]=\Mat{\Sigma}$ as before,    and specialize to mutually orthogonal steering vectors $ \Mat{A}_\mathcal{M}$  and Gauss loss,  i.e.,  $\Mat{R}'_{\Mat{Y}}=\Mat{0}$.  
This gives
\begin{align}
{\Vec{b}_m^{\sf H}  \mathsf{E}[\Mat{T}_{\Mat{Y}}] \Vec{b}_m} \gamma_m 
                   &= \Vec{a}_m^{\sf H} \Mat{\Sigma}^{-1}\Vec{a}_m   \gamma_m 
\end{align}
with
\begin{align}
 \Mat{\Sigma}&= \Mat{A}_\mathcal{M} \Mat{\Gamma}_\mathcal{M}  \Mat{A}_\mathcal{M}^{\sf H} + \sigma^2 \Mat{I}_N,
 \\
 \Mat{\Sigma}^{-1}&=\Mat{A}_\mathcal{M} \Mat{\Xi}_\mathcal{M} \Mat{A}_\mathcal{M}^{\sf H} + \frac{1}{\sigma^2} \Mat{I}_N
 \\
\Mat{\Xi}_\mathcal{M} &= \frac{1}{N}\mathop{\mathrm{diag}}\left(\frac{1}{\sigma^2+N\gamma_m} - \frac{1}{\sigma^2}\right)_{m\in\mathcal{M}}
\end{align}
and
\begin{align}
          \mathsf{E}[\beta]         &= 
\max_{m\in\mathcal{M}} \left(  \Vec{a}_m^{\sf H} \Mat{\Sigma}^{-1}\Vec{a}_m   \gamma_m \right)   \nonumber\\ 
&= \max_{m\in\mathcal{M}} \left(\frac{N\gamma_m}{\sigma^2+N\gamma_m} \right).
\end{align}
For $ \max_m (N\gamma_m) \ge \sigma^2$,  i.e.  at ASNR $\ge0\,$dB, we see that 
$\frac12\le\mathsf{E}[\beta]<1$.

\subsection{CRB for multiple DOA estimation}
\label{sec:CRB}
\subsubsection{Gaussian array data model}
\label{sec:Gaussian-CRB}
The CRB for multiple DOA estimation from Gaussian array data $\Mat{Y}$ is evaluated according to \cite[Eq.~(8.106)]{VanTreesBook}.
In Figures showing RMSE vs. ASNR, the trace of \eqref{eq:Gaussian-CRB} is plotted for $\sigma^2 = \frac{N}{\mathrm{ASNR}}$,

\newcommand{\GammaM}{\Mat{\Gamma}_{\mathcal{M}}}
\newcommand{\AM}{\Mat{A}_{\mathcal{M}}}
\newcommand{\RealBig}[1]{\mathop{\mathrm{Re}}\Big\{ #1 \Big\}}

\begin{align}
  C_{\mathrm{CR,Gauss}}(\Vec{\theta}) = \frac{\sigma^2}{2L} \mathop{\mathrm{tr}} & \Bigg\{ \RealBig{ \Big[  \GammaM 
                \Big( \Mat{I}_K + \AM^{\sf H} \AM \frac{\GammaM}{\sigma^2} \Big)^{-1} \nonumber \\
          &      \Big( \AM^{\sf H} \AM \frac{\GammaM}{\sigma^2} \Big)
            \Big] \odot \Mat{H}^{\sf T}}^{-1} \Bigg\} \label{eq:Gaussian-CRB}
\end{align}
with 
\begin{align}
    \GammaM &= \mathop{\mathrm{diag}}(\Vec{\gamma}_{\mathcal{M}}), \\
    \Mat{H} &= \Mat{D}^{\sf H} \left( \Mat{I}_N - \AM \AM^+ \right) \Mat{D}, \\
    \Mat{D} &= \left[ \left( \frac{\partial \Vec{a}(\theta)}{\partial \theta} \right)_{\theta=\theta_1}
    \ldots \hspace{2ex}
\left( \frac{\partial \Vec{a}(\theta)}{\partial \theta} \right)_{\theta=\theta_K} \right].
\end{align}

\subsubsection{CES array data model}
\label{sec:CES-CRB}

The CRB for multiple DOA estimation from CES distributed array data  $\Mat{Y}$ is derived in \cite[Eq. (20)]{besson2013fisher} and \cite[Eq. (17)]{greco2013cramer} 
based on the Slepian-Bangs formula.  
Starting from $p_{\Vec{y}}(\Vec{y}) = p_{\Vec{y}}(\Vec{y} | \Vec{\theta})$ given the true source scenario $\Vec{\theta}$ as defined in \eqref{eq:ces},  this gives 
\begin{align}
  C_{\mathrm{CR,CES}}(\Vec{\theta}) &= \frac{1}{L} \mathop{\mathrm{tr}} \{ \Mat{F}^{-1}\},
\label{eq:CES-CRB}
\end{align}
where the Fisher information matrix $\Mat{F}$ has elements
\begin{align} 
F_{i,j} = & \mathsf{E} \left[  \frac{ \partial \log p_{\Vec{y}}(\Vec{y}| \Vec{\theta})}{ \partial \theta_i}  \frac{ \partial \log p_{\Vec{y}}(\Vec{y}| \Vec{\theta})}{ \partial \theta_j}   \right]   \\
    = &( \psi_1 - 1 ) \mathrm{tr}(\Mat{\Sigma}^{-1} \Mat{\Sigma}_i) \; \mathrm{tr}(\Mat{\Sigma}^{-1} \Mat{\Sigma}_j)   
 \quad +  \nonumber \\
& \hspace{5.5ex} \psi_1  \mathrm{tr}(\Mat{\Sigma}^{-1} \Mat{\Sigma}_i   \Mat{\Sigma}^{-1} \Mat{\Sigma}_j) 
 \label{eq:FIM-CES}
\end{align} 
with array covariance matrix $\Mat{\Sigma}$ defined in \eqref{eq:Sigma-model_new},  $\Mat{\Sigma}_i = \frac{\partial \Mat{\Sigma}}{\theta_i}$ and 
\begin{equation} \label{eq:psi1}
\psi_1 =   \frac{\mathsf{E}[ \psi(\Vec{y}^{\sf H} \Mat{\Sigma}^{-1} \Vec{y} )^2]}{N(N+1)}.
\end{equation} 
For the MVT array data model,  this evaluates to
\begin{equation} \label{eq:psi1MVT}
\psi_1^{\mathrm MVT} =   \frac{2N + \nu_{\mathrm{data}}}{2(N+1)+ \nu_{\mathrm{data}}}
\end{equation} 
and the MVT CRB, $C_{\mathrm{CR,MVT}}(\Vec{\theta})$,  is evaluated by \eqref{eq:CES-CRB} with \eqref{eq:FIM-CES} and \eqref{eq:psi1MVT}.

\end{document}